\documentclass[a4paper,fleqn]{cas-sc}

\usepackage[authoryear,longnamesfirst]{natbib}
\usepackage{orcidlink} 

\usepackage{xcolor}
\usepackage{graphicx}
\usepackage{booktabs}
\usepackage{cancel}
\usepackage{algorithmic}
\usepackage{pgfplotstable}
\usepackage[nameinlink,noabbrev]{cleveref}
\usetikzlibrary{external}
\usepackage{tikz}
\usepackage{pgfplots}
\usetikzlibrary{shapes.arrows, patterns, calc}
\usepackage{tikz-3dplot}
\ifpdf
  \DeclareGraphicsExtensions{.eps,.pdf,.png,.jpg}
\else
  \DeclareGraphicsExtensions{.eps}
\fi

\usepgfplotslibrary{groupplots}
\usetikzlibrary{arrows}
\usetikzlibrary{shapes}
\usetikzlibrary{decorations.text}
\usetikzlibrary{spy}
\usetikzlibrary{calc}
\usetikzlibrary{positioning}
\usepackage{siunitx}
\usepackage{comment}
\usetikzlibrary{arrows.meta}
\definecolor{steelblue}{HTML}{A1BDC7}
\definecolor{orange}{HTML}{D98C21}
\definecolor{silver}{HTML}{B0ABA8}
\definecolor{rust}{HTML}{B8420F}
\definecolor{seagreen}{HTML}{2E6B69}
\definecolor{joshua}{HTML}{FBDC7F}
\definecolor{darksky}{HTML}{154c79}

\colorlet{lightsilver}{silver!30!white}
\colorlet{lightlightsilver}{silver!10!white}
\colorlet{darkorange}{orange!85!black}
\colorlet{darksilver}{silver!85!black}
\colorlet{darksteelblue}{steelblue!85!black}
\colorlet{darkrust}{rust!85!black}
\colorlet{darkseagreen}{seagreen!85!black}
\colorlet{darkjoshua}{joshua!85!black}

\definecolor{gblue}{HTML}{1D428A}
\definecolor{bgreen}{HTML}{00471B} 
\definecolor{hred}{HTML}{C8102E}  
\definecolor{cgreen}{HTML}{007A33}  

\definecolor{fred}{HTML}{B40F20}
\definecolor{fdarkorange}{HTML}{E58606}

\definecolor{ztealblue}{HTML}{3B9AB2}   
\definecolor{zlightblue}{HTML}{78B7C5}  
\definecolor{zyellow}{HTML}{EBCC2A}     
\definecolor{zgolden}{HTML}{E1AF00}     
\definecolor{zred}{HTML}{F21A00}

\colorlet{DGColor}{joshua}
\colorlet{LimiterColor}{orange}
\colorlet{IGRrTwoColor}{rust}
\colorlet{IGRrFiveColor}{seagreen}
\colorlet{ExactColor}{steelblue}
\colorlet{RefColor}{steelblue}

\definecolor{reflightgray}{HTML}{8C8C8C} % light charcoal
\definecolor{reflightblue}{HTML}{8FAFC4} % light blue–gray
\definecolor{reflightgreen}{HTML}{A6BE8A} % light olive
\definecolor{reflightorange}{HTML}{D1A089} % light brick

\colorlet{RefOneColor}{reflightgray}    % m^{-1}
\colorlet{RefTwoColor}{reflightblue}    % m^{-2}
\colorlet{RefThreeColor}{reflightorange} % m^{-3}
\colorlet{RefFourColor}{reflightgreen} % m^{-4}

\colorlet{p1color}{rust} % p = 1
\colorlet{p3color}{seagreen} % p = 3
\colorlet{p5color}{darksky} % p = 5
\colorlet{p7color}{darksky} % p = 7

\pgfplotsset{ compat=1.18,
    standard/.style={
      scale only axis,
      width=0.15\textwidth,
      height=0.15\textwidth,
      enlarge x limits=0.05,
      enlarge y limits=0.05,
      max space between ticks=40,
      every axis/.append style={font=\small},
      every legend/.append style={font=\small},
      every node/.append style={font=\small},    
      },
      3d/.style={
      colormap={darkskycolormap}{
        color(0.0)=(zred),  
        color(0.05)=(zyellow),
        color(0.5)=(cgreen),
        color(1.0)=(gblue),
      }
      },
}

\newcommand{\pgfplotscomputeslope}[3]{%
  \pgfplotstableread[col sep=comma]{#1}\pgfplotslocaltable
  \pgfplotstablecreatecol[
    create col/linear regression={
      table=\pgfplotslocaltable,
      x=m,
      y=#2,
      xmode=log,
      ymode=log,
    }
  ]{regression_column}{\pgfplotslocaltable}%
  \expandafter\edef\csname #3\endcsname{\pgfplotstableregressiona}%
}

\pgfplotsset{
colormap={inferno}{
rgb255=(0,0,4)
rgb255=(31,12,72)
rgb255=(85,15,109)
rgb255=(136,34,106)
rgb255=(186,54,85)
rgb255=(227,89,51)
rgb255=(249,140,10)
rgb255=(249,201,50)
rgb255=(252,255,164)
}
}
\usepackage{mathtools}
\usepackage{amssymb}
\usepackage{bm}

\mathtoolsset{centercolon}

\renewcommand{\phi}{\varphi}
\renewcommand{\epsilon}{\varepsilon}

\newcommand{\Id}{\mathbf{I}}

\newcommand{\vct}[1]{\bm{#1}}

\newcommand{\trace}{\operatorname{tr}}

\DeclarePairedDelimiterX{\infdivx}[2]{(}{)}{#1\;\delimsize\|\;#2}

\DeclarePairedDelimiterX{\Biginfdivx}[2]{\Big(}{\Big)}{#1\;\delimsize\big\|\;#2}

\DeclarePairedDelimiterX{\biginfdivx}[2]{\big(}{\big)}{#1\;\delimsize\|\;#2}

\newcommand\Ccancel[2][black]{%
    \let\OldcancelColor\CancelColor
    \renewcommand\CancelColor{\color{#1}}%
    \cancel{#2}%
    \renewcommand\CancelColor{\OldcancelColor}%
}

\ifpdf
  \DeclareGraphicsExtensions{.eps,.pdf,.png,.jpg}
\else
  \DeclareGraphicsExtensions{.eps}
\fi

\newtheorem{example}{Numerical Example}

\crefname{hypothesis}{Hypothesis}{Hypotheses}
\crefname{example}{Example}{Examples}

\begin{document}

\let\WriteBookmarks\relax
\def\floatpagepagefraction{1}
\def\textpagefraction{.001}

\shorttitle{IGR--DG Semidiscretization of Euler equations}
\shortauthors{Eyob, Arias, Bryngelson, and Sch{\"a}fer}

\title[mode=title]{Discontinuous Galerkin Semidiscretization of the Information Geometric Regularized Compressible Euler Equations}

\author[1]{Brook Eyob}
\ead{brook@gatech.edu}

\author[1]{Jesus Arias}
\ead{jarias9@gatech.edu}

\author[1]{Spencer H.\ Bryngelson}
\ead{shb@gatech.edu}

\author[2]{Florian Sch{\"a}fer}
\ead{florian.schaefer@nyu.edu}

\affiliation[1]{
  organization={Georgia Institute of Technology},
  city={Atlanta},
  state={Georgia},
  country={USA}
}

\affiliation[2]{
  organization={New York University},
  city={New York},
  state={New York},
  country={USA}
}

\begin{abstract}
Shock stabilization in compressible Euler flows remains a central challenge for high-order numerical methods.
Existing shock-capturing approaches, including limiters, artificial viscosity, and reconstruction-based methods, involve tradeoffs between robustness, accuracy, preservation of fine-scale flow features, and computational complexity.
In this work, we develop a discontinuous Galerkin (DG) discretization of the information geometric regularization (IGR) framework introduced by Cao and Sch{\"a}fer for the compressible Euler equations.
The method regularizes shocks at the PDE level by augmenting the Euler equations with the entropic pressure $\Sigma$, obtained from an auxiliary elliptic equation.
Within the DG formulation, the regularization enters only through the augmented pressure $P+\Sigma$ in the Euler fluxes, preserving the conservative structure of the discretization while using a common approximation space for both the hyperbolic and elliptic equations.
Numerical experiments spanning one and two-dimensional benchmark problems show the proposed formulation stabilizes shocks without shock-capturing limiters or artificial viscosity, although positivity-preserving methods may still be required when the density or pressure approaches zero.
Compared with a characteristic TVB-limited DG formulation, the IGR--DG method resolves increasingly finer-scale flow features as the polynomial order is increased while maintaining stable shock resolution.
The entropic pressure remains localized to regions of strong compression with minimal activation in smooth regions of the flow, providing selective PDE-level regularization while preserving the underlying solution elsewhere.
\end{abstract}


% 65M60, 35L65, 76M10, 65M12, 35Q31, 35L67
% 65M60 - Galerkin methods
% 35L65 - Conservation laws
% 76M10 - Finite element Fluid Mechanics 
% 65M12 - Stability and convergence
% 35Q31 - Euler Equations
% 35L67 - Shocks and singularities for hyperbolic equations

\maketitle

\section{Introduction}
The compressible Euler equations describe the conservation of mass, momentum, and energy in compressible fluids of negligible viscosity.
They arise in applications ranging from high-speed aerodynamics~\citep{anderson2002modern} and atmospheric dynamics~\citep{vallis2017atmospheric} to astrophysical flows~\citep{shu1991physics,pringle2007astrophysical}.
The nonlinear hyperbolic structure of the Euler equations supports wave propagation and complex flow interactions.

A defining feature of these equations is the formation of shock waves from smooth initial conditions, and numerical methods must resolve both smooth flow features and shocks within a single discretization.
Developing methods that balance these requirements remains a central challenge in computational fluid dynamics.

\subsubsection*{Shock formation}

Nonlinear effects can cause smooth solutions of the compressible Euler equations to develop singularities in finite time.
In compressive regions, characteristics converge as faster-moving waves overtake slower-moving waves.
This steepening produces unbounded solution gradients, causing the classical solution to break down and shocks to form~\citep{lax1964development,smoller1994shock,buckmaster2023formation}.
Beyond this point, the equations admit weak solutions that may contain jump discontinuities~\citep{dafermos2005hyperbolic}.

\subsubsection*{Existing shock-capturing strategies}

A wide range of numerical techniques exist to stabilize shocks in the compressible Euler equations, including artificial viscosity, limiter-based methods, adaptive reconstruction, filtering, and structure-preserving formulations.

Artificial viscosity methods introduce localized diffusion near steep gradients to suppress spurious oscillations.
Classical approaches relied on grid-scale features~\citep{vonneumann1950method,jameson1981numerical,cook2005hyperviscosity,mani2009suitability}, while modern techniques use more sophisticated indicators such as subcell modal energy decay~\citep{persson2006sub}, auxiliary regularization equations~\citep{barter2010shock}, or monitor local violations of the entropy inequality~\citep{guermond2011entropy,chan2025artificial}.

Limiter-based approaches suppress numerical oscillations by clamping or smoothing out steep gradients near shocks.
These strategies include MUSCL and TVD schemes~\citep{van1979muscl, harten1983high} and generalized slope and flux limiters designed for high-order methods~\citep{cockburn1998runge, kuzmin2010vertex,hennemann2020entropy}.

ENO and WENO methods adapt reconstruction stencils using local smoothness indicators to improve resolution near discontinuities while avoiding spurious oscillations~\citep{harten1987eno,jiang1996weno}.

Positivity preserving and invariant-domain preserving methods enforce physical admissibility by maintaining nonnegative density and pressure and keeping the solution within an admissible state set~\citep{zhang2010positivity,guermond2018second,kuzmin2020monolithic,pazner2021sparse}.

Filtering and spectral viscosity methods damp high-frequency modes while preserving larger-scale flow structure~\citep{tadmor1989ssv,hesthaven2007nodal}.
Entropy-stable and entropy-consistent formulations impose additional thermodynamic constraints on the discrete solution~\citep{tadmor1987entropy,ismail2009entropy,chandrashekar2013kinetic}.
Finally, some approaches aim to track the geometric structure of the shock and align it with DG cell interfaces \cite{zahr2019implicit,naudet2024space,huang2022robust,anton2026bow}.

\subsubsection*{Limitations of existing shock-capturing strategies}

Existing shock-capturing approaches provide robust and accurate simulations of compressible flows, but each introduces tradeoffs among accuracy, robustness, and computational complexity~\citep{hoskin2024discontinuous}.

Artificial viscosity methods regularize discontinuities by introducing additional dissipation.
The amount and localization of this dissipation must balance shock resolution against the preservation of fine-scale features.
In high-Mach-number regimes, many methods based on artificial viscosity suffer from a lack of robustness \cite{peck2025comparing}.

Limiter-based approaches can preserve accuracy under suitable parameter choices~\citep{cockburn2001runge}, in practice however, they may activate near smooth and physically relevant structures and damp fine-scale features, especially in high-order methods~\citep{trojak2021shock}.
This shortcoming is partially addressed by subcell limiting strategies~\citep{dumbser2014posteriori,hennemann2020entropy}.

ENO and WENO methods improve shock resolution by adaptively reconstructing stencils using local smoothness indicators.
However, the reconstruction procedure adds algorithmic complexity and becomes more difficult to implement on nonuniform, unstructured, or multidimensional meshes~\citep{shu2009high, wang2013highorder}.

Filtering and spectral-viscosity methods selectively damp high-frequency modes to improve stability.
Although they can be less intrusive than conventional artificial viscosity, they still introduce dissipation and may attenuate physically relevant structures when applied aggressively~\citep{xu2018spectral}.

Positivity-preserving and invariant-domain-preserving methods maintain physically admissible discrete states and improve robustness under challenging flow conditions.
These methods primarily enforce admissibility and often require additional shock-capturing mechanisms, since invariant-domain preservation alone does not prevent Gibbs oscillations in higher-order schemes~\citep{zhang2011maximum,guermond2016invariant}.

Entropy-stable and entropy-consistent formulations enforce discrete thermodynamic constraints and improve nonlinear robustness.
Entropy stability alone, however, does not remove the need for additional shock stabilization near under-resolved discontinuities~\citep{hiltebrand2014entropy,hennemann2020entropy}.

Excplicitly or implicitly, all of the above approaches rely on a form of artificial diffusion and therefore face a fundamental trade-off between stability and accuracy.
Shock tracking methods such as \cite{zahr2019implicit} can avoid this trade-off, but are limited to sufficiently simple shock geometries.
This motivates the search for PDE-level regularization strategies that use mechanisms other than viscosity.
\subsubsection*{Information geometric regularization (IGR)}

Information geometric regularization (IGR), introduced by Cao and Sch{\"a}fer~\citep{cao_schafer2023barotropic}, is the first \emph{inviscid} regularization for Euler-like equations.
In the original formulation, derived for the barotropic case, the pressure in the momentum flux is augmented according to $P \mapsto P+\Sigma$, where the entropic pressure $\Sigma$ satisfies an auxiliary elliptic equation.
The regularization responds to compressive deformation while preserving the conservative structure of the governing equations.

In this work, we apply the augmented pressure $P+\Sigma$ to the full compressible Euler equations with an independent energy equation.
Relative to the barotropic formulation, this extension introduces an energy conservation law and replaces the classical energy flux $(E+P)\vct{u}$ with $(E+P+\Sigma)\vct{u}$.

In Eulerian variables, the IGR-regularized Euler equations take the form
\begin{equation}
\label{eq:igr-euler}
\begin{aligned}
&\partial_t
\begin{pmatrix}
\rho \\
\rho\vct{u} \\
E
\end{pmatrix}
+
\nabla\cdot
\begin{pmatrix}
\rho\vct{u} \\
\rho\vct{u}\otimes\vct{u}+(P+\Sigma)\Id \\
(E+P+\Sigma)\vct{u}
\end{pmatrix}
=
\begin{pmatrix}
0 \\
\vct{f} \\
\vct{f}\cdot\vct{u}
\end{pmatrix},
\\[0.75em]
&\rho^{-1}\Sigma
-
\alpha\,\nabla\cdot\left(\rho^{-1}\nabla\Sigma\right)
=
\alpha\left(\trace(D\vct{u})^2+\trace\left((D\vct{u})^2\right)\right).
\end{aligned}
\end{equation}

Here, $\rho$ denotes the fluid density, $\vct{u}$ the velocity field, $E$ the total energy density, $P$ the thermodynamic pressure, $\Sigma$ the entropic pressure, and $\vct{f}$ an external body force.
The parameter $\alpha>0$ controls the strength of the regularization and the characteristic length scale over which it redistributes compression.
The resulting equations replace shocks with smooth profiles and avoid the explicit viscous dissipation used by many classical regularization methods~\citep{cao_schafer2023barotropic}.
\cite{cao_schafer2024pressureless} proved rigorously that the IGR system admits global smooth solutions in the unidimensional pressureless case. 
\cite{barham2026shock} prove the existence of smooth traveling wave solutions of unidimensional IGR, for a large class of equations of state.

\cite{wilfong2025simulating,radhakrishnan2026shocks} apply IGR to finite volume methods, achieving the largest-ever compressible flow simulations at one quadrillion degrees of freedom.
Subsequent work has extended the IGR framework through Hamiltonian, thermodynamic, and directional formulations~\citep{taylor2026thermodynamically,xu2026compression}.

\subsubsection*{IGR for discontinuous Galerkin methods}

Discontinuous Galerkin methods provide a high-order framework for approximating hyperbolic conservation laws such as the compressible Euler equations~\citep{cockburn2001runge}.
By combining elementwise polynomial approximations with conservative numerical fluxes, DG methods accommodate jumps between elements while retaining high-order accuracy in smooth regions~\citep{hesthaven2007nodal}.
Shocks nevertheless generate spurious oscillations within the polynomial representation, so practical DG formulations require an additional stabilization mechanism~\citep{gottlieb1997gibbs}.

IGR is well suited to the DG setting because it preserves the conservative form of the Euler equations and enters the hyperbolic system only through the augmented pressure $P+\Sigma$.
Consequently, the existing DG transport discretization can be retained, while an auxiliary elliptic solve provides the stabilizing entropic pressure used in the momentum and energy fluxes.
This PDE-level regularization can stabilize shocks without suppressing the higher-order modes and fine-scale features that motivate the use of DG methods.

Earlier work illustrated IGR using finite-difference discretizations, and subsequent studies applied the method in finite-volume simulations of compressible flows~\citep{wilfong2025simulating,radhakrishnan2026shocks}.
The present work develops the first discontinuous Galerkin realization of IGR for the compressible Euler equations and examines how the regularization interacts with higher-order polynomial approximation and discontinuous interelement coupling.

We evaluate the resulting IGR--DG method on benchmark problems involving shock formation, shock propagation, shock--shock interaction, and shock--vortex interaction.
These experiments assess stability under compression, localization of the entropic pressure, and preservation of fine-scale structures, with a characteristic TVB-limited DG formulation serving as the baseline comparison.
We also examine how the required regularization strength varies with shock intensity, which provides practical guidance for selecting the regularization parameter.

We provide a Julia implementation as a compact research code for structured Cartesian meshes and an implementation in MFEM~\citep{mfem-2021,mfem-2024} that extends the formulation to curved unstructured meshes, more general geometries, and a broader range of boundary conditions.

\section{Information Geometric Regularization (IGR)}
Information geometric regularization (IGR) was introduced by Cao and Sch{\"a}fer as an inviscid regularization of the barotropic Euler equations based on ideas from information geometry, geometric hydrodynamics, and interior-point methods~\citep{cao_schafer2023barotropic,amari2016information,ebin1970groups,guler1996barrier,khesin2021geometric}.
The original work developed the geometric framework through the one-dimensional pressureless Euler equations before extending the construction to multidimensional barotropic flows~\citep{cao_schafer2023barotropic}.
Subsequent work established rigorous analytical results for the pressureless setting and further investigated the structure of the regularized dynamics~\citep{cao_schafer2024pressureless}.

\subsubsection*{Geometric perspective on shock formation}

The flow-map description of fluid motion provides a geometric interpretation of shock formation and motivates the IGR framework.
The flow map $\Phi_t$ carries each fluid particle from its initial position to its location at time $t$~\citep{arnold1966geometrie,ebin1970groups}.
In the pressureless case, the $\Phi_t$ follows a geodesic on the manifold of diffeomorphisms~\citep{khesin2021geometric}.
As long as the solution remains smooth, $\Phi_t$ remains invertible and distinct particle trajectories remain separated throughout the evolution.
The associated deformation gradient $D\Phi_t$ describes the local deformation of the flow. 
The Jacobian determinant $\det(D\Phi_t)$ measures the corresponding local volume change.

Compression drives the Jacobian determinant toward zero as neighboring fluid particles converge.
The loss of invertibility corresponds to particle trajectories intersecting.
At this point, the classical solution breaks down and a shock forms~\citep{khesin2007shock,lax1978accuracy,cao_schafer2024pressureless}.
This geometric interpretation motivates regularization mechanisms that prevent compressive deformations from reaching the singular limit.

\subsubsection*{IGR formulation}

Information geometric regularization prevents shock formation by penalizing compressive deformations before the flow map loses invertibility~\citep{cao_schafer2023barotropic}.
In the multidimensional setting, the resulting regularization takes the form of a logarithmic barrier functional,
\begin{equation}
    \psi_\alpha[\Phi,\Phi']
    =
    \frac{1}{2}\int_\Omega \|\Phi(x,t)-x\|^2\,dx
    +
    \alpha\int_\Omega -\log\!\left(\Phi'(x,t)\right)\,dx,
\end{equation}
where $\alpha>0$ controls the strength of the regularization and $\Phi'$ is the (formally independent) Jacobian determinant $\det(D\Phi)$.
As $\det(D\Phi)\to0$, the logarithmic term diverges and assigns an increasingly large cost to compression, preventing the deformation gradient from becoming singular.

Mass conservation provides an Eulerian interpretation of the barrier,
\begin{equation}
    \rho(\Phi(x,t),t)=\left(\det D\Phi(x,t)\right)^{-1},
\end{equation}
which relates the local compression of the flow map to the density field.
Substituting this relation into the barrier functional and performing a change of variables yields
\begin{equation}
    \int_\Omega -\log\!\left(\det D\Phi\right)\,dx
    =
    \int_\Omega \rho(x,t)\log\rho(x,t)\,dx,
\end{equation}
showing that the regularization penalizes compressive mass concentration through the Shannon entropy of the density.
Information geometric regularization replaces the $L^2$ geodesics defining the nominal equation with the dual geodesics defined by the barrier functional, in the sense of information geometry \cite{amari2016information}.

Expressing the resulting equations in Eulerian variables allows summarizing these geometric effects through the entropic pressure $\Sigma$.
Varying the functional produces the $\Sigma$, yielding the IGR regularization introduced in \Cref{eq:igr-euler}.
The entropic pressure is determined by the elliptic equation
\begin{equation}
\rho^{-1}\Sigma
-
\alpha\,\nabla\cdot\left(\rho^{-1}\nabla\Sigma\right)
=
\alpha\left(
\left(\trace(D\vct{u})\right)^2
+
\trace\left((D\vct{u})^2\right)
\right).
\label{eq:entropic-pressure}
\end{equation}
The left-hand side defines an elliptic operator for the entropic pressure field, while $\alpha$ controls both the strength of the regularization and the characteristic length scale over which deformation is redistributed.

The right-hand side depends on invariants of the velocity gradient tensor associated with volumetric and local flow deformation.
The first invariant, $\left(\trace(D\vct{u})\right)^2$, measures the magnitude of the local volumetric deformation, and always yields a positive contribution.
The second invariant, $\trace\left((D\vct{u})^2\right)$, captures additional local deformation, including strain and rotational effects, and is not restricted to positive values.
Together, these invariants allow the entropic pressure to respond to the complete local deformation of the flow.
The entropic pressure remains predominantly positive, with its strongest response occurring in regions dominated by compressive strain.

Since these invariants require a well-defined velocity gradient tensor, only sufficiently smooth initial conditions are admissible for the IGR formulation.
Discontinuous initial conditions are therefore smoothed before time integration so that they remain consistent with the IGR formulation.
The entropic pressure is computed directly from the instantaneous flow state rather than from an equation of state.
Although the entropic pressure enters the equations the like thermodynamic pressure $P$, it reflects geometric compression rather than thermodynamic effects.

\subsubsection*{Relation to existing regularization methods}

IGR differs from conventional shock-capturing methods by using an \emph{inviscid} mechanism. 
Most shock-capturing strategies use dissipative fluxes, limiting procedures, reconstruction-based corrections, or localized artificial viscosity to control oscillations near shocks~\citep{hesthaven2018numerical}.
IGR instead modifies the compressive dynamics through the entropic pressure while preserving the conservative form of the Euler equations~\citep{cao_schafer2023barotropic}.
As discussed by \cite[Section 7.1]{cao_schafer2024pressureless} it has the character of a nonlocal and sign-indefinite regularization.

This distinction places IGR among continuum regularization methods, although its mechanism differs from viscous and entropy-viscosity approaches that stabilize shocks through diffusion.
The resulting dynamics replace discontinuous shocks with smooth profiles whose characteristic width is proportional to $\sqrt{\alpha}$, without dissipating fine-scale features \citep{cao_schafer2023barotropic}.

\subsection{Properties of IGR regularization}

The regularized equations possess several structural properties that distinguish them from the classical Euler equations.
The following subsections describe the coupled hyperbolic--elliptic structure, the role of the regularization parameter, and the selective activation of the entropic pressure.

\subsubsection*{Hyperbolic--elliptic structure}

The IGR-regularized equations couple the hyperbolic Euler equations to the elliptic equation that determines the entropic pressure $\Sigma$~\citep{cao_schafer2024pressureless}.
The conservative variables $(\rho,\rho\vct{u},E)$ evolve through hyperbolic transport with modified fluxes, while the elliptic equation recovers $\Sigma$ from the instantaneous flow state.

The entropic pressure $\Sigma$ has no separate evolution equation, instead $\Sigma$ is determined from the flow state through the elliptic relation \eqref{eq:entropic-pressure}.
Because the regularization enters only through the augmented pressure $(P + \Sigma)$, the conservative transport structure of the Euler equations is preserved.

\subsubsection*{Regularization parameter $\alpha$}

The regularization parameter $\alpha$ controls both the strength and spatial extent of the entropic pressure response.
Larger values produce broader regularized shock profiles and distribute deformation over a greater distance, while smaller values localize the regularization more strongly near developing shocks.

An important property of the regularized equations is consistency with the classical Euler equations.
As $\alpha\to0$, the entropic pressure contribution vanishes  and the regularized equations approach the unregularized compressible Euler equations.
The corresponding regularized shock profiles narrow as the regularization weakens and approach the discontinuous shock structure associated with the Euler solution.

\subsubsection*{Selective $\Sigma$ activation}

The entropic pressure responds most strongly in regions of large deformation while remaining weak where the flow varies smoothly.
When the velocity gradients are small, the right-hand side of \cref{eq:entropic-pressure} is small and the dynamics remain close to those of the unregularized Euler equations.
As compression intensifies, the entropic pressure redistributes the deformation over a finite length scale and replaces the developing discontinuity with a narrow smooth profile.
The entropic pressure is determined directly by the elliptic equation and does not require an artificial shock sensor, limiting procedure, or auxiliary activation criterion.

\section{Information Geometrically Regularized Discontinuous Galerkin}
Discontinuous Galerkin (DG) methods provide a high-order discretization of the IGR-regularized equations while preserving their conservative transport structure.
Originally developed as conservative high-order discretizations for hyperbolic conservation laws, DG methods combine element-local polynomial approximations with weak interelement coupling~\citep{cockburn2001runge, hesthaven2007nodal}.

The IGR--DG method couples a DG discretization of the Euler equations with a compatible discretization of the elliptic equation for the entropic pressure~\citep{arnold2002unified,ern2006discontinuous}.
Both equations are discretized in the same discontinuous function space, providing a consistent hyperbolic--elliptic discretization.

\subsection{Euler equation DG discretization}

The DG discretization of the hyperbolic subsystem begins with the IGR-regularized equations written in conservative form as
\begin{equation}
\label{eq:igr_euler_conservative}
\partial_t U
+
\nabla \cdot F(U,\Sigma)
=
0,
\end{equation}
where the conserved variables are $U=(\rho,\rho \vct{u},E)^T$ and the IGR-regularized Euler flux is
\begin{equation}
F(U,\Sigma)
=
\begin{pmatrix}
\rho \vct{u} \\
\rho \vct{u} \otimes \vct{u} + (P+\Sigma) \Id \\
(E+P+\Sigma)\vct{u}
\end{pmatrix}.
\end{equation}

Let $\Omega\subset\mathbb{R}^d$ denote the computational domain partitioned into nonoverlapping elements $\{\Omega_k\}_{k=1}^{N_e}$.
On each element $\Omega_k$, the discrete solution is represented in a local basis $\{\phi_{k,j}\}_{j=1}^{N_k}$ as
\begin{equation}
\label{eq:dg_solution_expansion}
U_h(x,t)\big|_{\Omega_k}
=
\sum_{j=1}^{N_k}
U_{k,j}(t)\phi_{k,j}(x),
\end{equation}
where $U_{k,j}(t)$ are the DG coefficients associated with the conserved variables~\citep{hesthaven2007nodal}.

Let $V_h$ denote the corresponding discontinuous space of piecewise polynomial functions.
Multiplying \cref{eq:igr_euler_conservative} by an arbitrary test function $v\in V_h$ and integrating over each element results in the weak form of the hyperbolic equations,
\begin{equation}
\int_{\Omega_k}
\partial_t U_h\,v
\,dx
+
\int_{\Omega_k}
\nabla\cdot F(U_h,\Sigma_h)\,v
\,dx
=
0.
\end{equation}
Applying integration by parts elementwise yields
\begin{equation}
\label{eq:hyperbolic_dg_weak_form}
\int_{\Omega_k}
\partial_t U_h\,v \,dx
- \int_{\Omega_k} F(U_h,\Sigma_h)\cdot\nabla v \,dx
+ \int_{\partial\Omega_k} \left( F(U_h,\Sigma_h)\cdot n_k \right)v \,dx
= 0,
\end{equation}
where $n_k$ denotes the outward unit normal on $\partial\Omega_k$.

A numerical flux $\widehat{F}(U_h^-,U_h^+,\Sigma_h^-,\Sigma_h^+)$ resolves the interface contribution using the interior and exterior traces~\citep{toro2013riemann}, because the physical flux is not uniquely defined on element boundaries.
Using the same basis for the trial and test spaces, set $v=\phi_{k,i}$, and replace the physical flux with the numerical flux to obtain
\begin{equation}
\label{eq:hyperbolic_element_dg}
\int_{\Omega_k}
\partial_t U_h\,\phi_{k,i} \,dx
- \int_{\Omega_k} F(U_h,\Sigma_h)\cdot\nabla\phi_{k,i} \,dx
+ \int_{\partial\Omega_k} \widehat{F}(U_h^-,U_h^+,\Sigma_h^-,\Sigma_h^+) \cdot n_{k} \,\phi_{k,i} \,dx
= 0,
\end{equation}
for $k=1,\ldots,N_e$ and $i=1,\ldots,N_k$.

Substituting \cref{eq:dg_solution_expansion} into \cref{eq:hyperbolic_element_dg} yields
\begin{equation}
\label{eq:hyperbolic_semidiscrete}
M_k\partial_t U_k
- K_k(U_h,\Sigma_h)
+ B_k(U_h,\Sigma_h)
= 0,
\end{equation}
where $M_k$, $K_k$, and $B_k$ denote the elemental mass, stiffness, and boundary contributions, respectively.
Rearranging gives the system of ordinary differential equations
\begin{equation}
\label{eq:hyperbolic_semidiscrete_rhs}
\partial_t U_k 
= M_k^{-1}
\left(
K_k(U_h,\Sigma_h)
- B_k(U_h,\Sigma_h)
\right).
\end{equation}

\subsection{Elliptic equation DG discretization}

The elliptic equation in \cref{eq:entropic-pressure} is discretized using the same approximation space as the Euler equations.
On each element $\Omega_k$ the discrete entropic pressure is represented in the local basis as
\begin{equation}
\label{eq:elliptic_sigma_expansion}
\Sigma_h(x)\big|_{\Omega_k}
=
\sum_{j=1}^{N_k}
\sigma_{k,j}\phi_{k,j}(x),
\end{equation}
where $\sigma_{k,j}$ are the DG coefficients.

Multiplying \cref{eq:entropic-pressure} by a test function $v \in V_h$ and integrating over each element yields
\begin{equation}
\int_{\Omega_k} \rho^{-1}\Sigma_h\,v \,dx
- \alpha \int_{\Omega_k} \nabla\!\cdot\!\left(\rho^{-1}\nabla\Sigma_h\right) \,v \,dx
=
\alpha \int_{\Omega_k} \left( (\trace(D\vct{u}_h))^2+ \trace((D\vct{u}_h)^2)\right) \,v \,dx.
\end{equation}

Integrating the divergence term by parts yields
\begin{equation}
\label{eq:elliptic_element_weak}
\int_{\Omega_k}\rho^{-1}\Sigma_h\,v \,dx
+ \alpha \int_{\Omega_k} \rho^{-1}\nabla\Sigma_h\cdot\nabla v \,dx
- \alpha \int_{\partial\Omega_k} \left(\rho^{-1}\nabla\Sigma_h\cdot n_k\right)v \,dx
=
\alpha \int_{\Omega_k} \left( (\trace(D\vct{u}_h))^2 + \trace((D\vct{u}_h)^2) \right) v \,dx.
\end{equation}

The symmetric interior penalty Galerkin (SIPG) method resolves the interface contributions using a numerical  trace $\widehat{q}_{\Sigma,\phi_{k,i}}$  that account for the consistency, symmetric consistency, and penalty contributions at element interfaces~\citep{arnold2002unified}.

Using the same basis $\{ \phi_{k,i} \}_{i=1}^{N_k}$ for the trial and test spaces, the SIPG formulation yields
\begin{equation}
\label{eq:elliptic_element_dg}
\int_{\Omega_k} \rho^{-1}\Sigma_h\,\phi_{k,i} \,dx
+ \alpha \int_{\Omega_k} \rho^{-1}\nabla\Sigma_h\cdot\nabla\phi_{k,i} \,dx
- \alpha \int_{\partial\Omega_k} \widehat{q}_{\Sigma,\phi_{k,i}} \,dx
= \alpha \int_{\Omega_k} \left( (\trace(D\vct{u}_h))^2 + \trace((D\vct{u}_h)^2) \right)\phi_{k,i} \,dx,
\end{equation}
for $k=1,\ldots,N_e$ and $i=1,\ldots,N_k$.

Assembling the elemental contributions over all elements yields the global linear system
\begin{equation}
\label{eq:elliptic_global_system}
A(U)\sigma = b(U),
\end{equation}
where $A(U)$ is the matrix associated with the discretized elliptic operator and $b(U)$ contains the discretized right-hand side of the elliptic equation~\citep{riviere2008discontinuous}.

\subsection{Coupled IGR--DG formulation}

The IGR--DG method couples the DG discretization of the Euler equations with the SIPG discretization of the elliptic equation.
The resulting semidiscrete system is
\begin{align}
\label{eq:coupled_semidiscrete_hyperbolic}
\partial_t U
&=
M^{-1}
\left(
K(U,\Sigma)
-
B(U,\Sigma)
\right),
\\
\label{eq:coupled_semidiscrete_elliptic}
A(U)\sigma
&=
b(U).
\end{align}

The matrix $A(U)$ depends on the discrete density field, while $b(U)$ depends on velocity gradients extracted from the discrete solution.
The resulting entropic pressure enters the Euler fluxes through the modified pressure contribution $P+\Sigma$.
The Euler discretization couples neighboring elements through numerical fluxes, while the SIPG discretization yields a globally coupled linear system.
At each time step, the elliptic system is solved using the current discrete solution and the resulting entropic pressure is used to advance the Euler equations.

\section{Implementation}
An implementation of the IGR–DG formulation that combines the DG discretizations of the Euler and elliptic equations within a single numerical framework is presented here.
Within this framework, the Euler equations are discretized using a standard DG formulation~\citep{hesthaven2007nodal}, with the entropic pressure entering through the numerical fluxes.
The elliptic equation is discretized using a symmetric interior penalty Galerkin (SIPG) formulation~\citep{riviere2008discontinuous} to compute the discrete entropic pressure field.
The method advances the conserved variables through explicit time integration and solves the elliptic problem at each time step to update the entropic pressure, coupling the regularization to the evolving flow while preserving the DG transport structure.

\subsection{Spatial discretization}

The DG approximation is constructed on a uniform Cartesian mesh of tensor-product elements in $d$ dimensions.
The domain is partitioned into nonoverlapping cells $\{\Omega_k\}^{N_e}_{k=1}$, each of the form
\begin{equation}
    \Omega_k = \prod_{\ell=1}^d [x_{\ell,k}^-,\,x_{\ell,k}^+],
\end{equation}
with identical side lengths in every coordinate direction.
The uniform, axis-aligned structure of the mesh implies that geometric factors are constant across elements, allowing basis functions, quadrature nodes, and gradient operators to be defined directly in physical coordinates and reused on every element.
On each element $\Omega_k$, the local approximation space is spanned by tensor-product Lagrange basis functions associated with a set of interpolation nodes $\{x_j\}_{j=1}^{N_k}$.
These basis functions are given by
\begin{equation}
\phi_j(x_1,\dots,x_d)
=
\prod_{\ell=1}^d
\left(
\prod_{\substack{m=1 \\ m\neq j_\ell}}^{n}
\frac{x_\ell - x_{\ell,m}}{x_{\ell,j_\ell} - x_{\ell,m}}
\right),
\end{equation}
where $j=(j_1,\dots,j_d)$ indexes the multi-dimensional interpolation node and $n$ denotes the number of one-dimensional nodes per coordinate direction.
All volume and surface integrals appearing in the discrete formulation are evaluated using Gauss-Legendre quadrature.
The discrete representations of the conservative variables and the entropic pressure follow directly from the DG expansions introduced in the formulation.
The conservative variables are approximated according to \eqref{eq:dg_solution_expansion}, with $U=(\rho,\rho\vct{u},E)^\top$, while the entropic pressure is represented using the approximation in \eqref{eq:elliptic_sigma_expansion}.
The corresponding coefficients $U_{k,j}(t)$ and $\sigma_{k,j}(t)$ define the finite-dimensional degrees of freedom associated with the hyperbolic and elliptic equations, respectively.

\subsection{Numerical fluxes}

The DG approximation is discontinuous across element interfaces, requiring numerical fluxes to define unique interface states and provide conservative coupling between neighboring elements.
This work uses a local Lax--Friedrichs numerical flux for the hyperbolic Euler equations and SIPG numerical traces for the elliptic equation.
The average operator $\{\!\!\{\cdot\}\!\!\}$ is defined as the arithmetic mean of interior and exterior traces, and the jump operator $[\![\cdot]\!]$ denotes the difference of traces taken in the normal direction.

\subsubsection*{Euler Lax--Friedrichs flux}

For the Euler equations, the local Lax--Friedrichs numerical flux is
\begin{equation}
\label{eq:euler_lf_flux}
\widehat{F}
= \{\!\!\{ F(U,\Sigma) \}\!\!\}\cdot n
- \frac{\lambda}{2}\,[\![ U ]\!],
\end{equation}
where $U^{-}$ and $U^{+}$ denote the interior and exterior traces of the conservative variables and $\Sigma^{-}$, $\Sigma^{+}$ the corresponding traces of the entropic pressure.

The average flux term provides a consistent approximation of the physical Euler flux across the interface, while the jump term supplies numerical dissipation.
The parameter $\lambda$ is a local estimate of the maximum characteristic speed of the compressible Euler equations in the normal direction.
For an ideal gas, the characteristic speeds are approximated by $u_n \pm c$, where $u_n=\vct{u}\cdot n$ and $c=\sqrt{\gamma P/\rho}$ is the sound speed computed from the thermodynamic pressure.
Following the standard local Lax--Friedrichs formulation, the dissipation coefficient is chosen as $ \lambda = \left( |u_n^{-}| + c^{-}, \; |u_n^{+}| + c^{+} \right)$.
This construction yields a consistent and conservative numerical flux with dissipation determined by the local characteristic speeds.

\subsubsection*{Elliptic trace}

The SIPG discretization introduces numerical traces associated with the consistency and symmetric consistency terms at element interfaces.
On a face with unit normal $n$, these traces are defined as~\citep{riviere2008discontinuous}
\begin{equation}
\label{eq:elliptic_trace_definition}
\widehat{q}_{\Sigma_{h},v}
= \{\!\!\{\rho^{-1}\nabla\Sigma_h\cdot n\}\!\!\}[\![v]\!]
+ \{\!\!\{\rho^{-1}\nabla v\cdot n\}\!\!\}[\![\Sigma_h]\!]
- \eta[\![\Sigma_h]\!][\![v]\!].
\end{equation}
where $\eta$ is a penalty parameter that ensures stability of the discrete formulation, weakly enforces continuity of the discrete entropic pressure across element interfaces, and scales as $\eta \sim \frac{p^2}{h}$~\citep{epshteyn2007estimation}.

\subsection*{Equation of state and energy treatment}

The equation of state closes the compressible Euler equations by relating density, pressure, and energy.
Both isentropic and ideal-gas equations of state were implemented.

For the isentropic Euler equations, the pressure is given by
\begin{equation}
    P(\rho) = K \rho^\gamma,
\end{equation}
and only the mass and momentum equations are evolved.
The total energy is reconstructed from the density and velocity as
\begin{equation}
    E = \frac{P}{\gamma - 1} + \frac{1}{2}\rho |\mathbf{u}|^2.
\end{equation}
This formulation is computationally less intensive and is appropriate for smooth flows in which entropy remains constant.

For the ideal gas Euler equations, the full energy equation is evolved and the pressure is computed from the conserved variables,
\begin{equation}
    P = (\gamma - 1)\Big(E - \tfrac{1}{2}\rho |\mathbf{u}|^2\Big).
\end{equation}
This formulation accounts for shock heating and entropy production and is used for flows involving strong compression.

When IGR is applied, the pressure is replaced by $P + \Sigma$ in the momentum and energy fluxes, while the equation of state remains unchanged~\citep{cao_schafer2023barotropic}.
Unless otherwise specified, all experiments use $\gamma = 1.4$.

\subsection{Boundary conditions}

Boundary conditions enter the DG formulation through the numerical flux, which depends on the interior trace and the prescribed boundary treatment.
Each physical boundary condition is enforced by prescribing an exterior state or directly assigning the boundary numerical flux~\citep{hesthaven2007nodal}.

Periodic boundaries couple the numerical flux across paired faces, while Dirichlet boundaries impose a prescribed value through the ghost state.
Reflective boundaries impose zero mass and energy fluxes, with the augmented pressure contributing to the normal momentum flux, while solid no-slip walls eliminate the velocity at the wall.
Transmissive boundaries use characteristic information to distinguish inflow, subsonic outflow, and supersonic outflow, while zero-Neumann boundaries copy the numerical flux from the adjacent interior interface.

The elliptic problem uses periodic coupling across paired faces or a homogeneous Neumann boundary condition for $\Sigma$ at nonperiodic boundaries, enforced weakly by omitting the boundary interface contributions~\citep{riviere2008discontinuous}.

\begin{table}[pos=h!]
\centering
\begin{tabular}{|c|c|}
\hline
\textbf{BC Type} & \textbf{Boundary Treatment} \\ \hline
Periodic & copy from paired boundary \\ \hline
Dirichlet & $U^{+}=U_{\mathrm{bc}}$ \\ \hline
Reflective & $\widehat{F}=(0,P+\Sigma,0)^T$ \\ \hline
No-slip & $\vct{u}=0$ at the wall \\ \hline
Transmissive & characteristic inflow--outflow condition \\ \hline
Zero Neumann & copy adjacent interior numerical flux \\ \hline
Elliptic $\Sigma$ & $\nabla\Sigma\cdot n=0$ at nonperiodic boundaries \\ \hline
\end{tabular}
\caption{
Boundary condition implementation for the Euler and elliptic equations.
Depending on the boundary type, the method prescribes an exterior state, directly assigns the Euler boundary flux, or modifies the elliptic interface contributions.
}
\label{tab:bc_ghost_states}
\end{table}

\subsection{TVB limiter baseline}

A DG method with a TVB slope limiter (DG+limiter) provides a baseline comparison with the IGR--DG approach.
The implementation uses a characteristic-based TVB slope limiter~\citep{krivodonova2007limiters} in both one and two dimensions.
In all cases, the limiter acts locally on the discrete solution representation and does not modify the governing equations, in contrast to the regularization introduced by the IGR method.

The limiter operates on the DG solution within each element by reconstructing local solution slopes and projecting them onto the characteristic fields of the Euler equations evaluated at the cell-average state.
Let $U_h$ denote the DG approximation restricted to a given element with cell average $\bar U$.
A candidate DG slope $s_{\mathrm{DG}}$ is reconstructed from the element degrees of freedom using precomputed basis-dependent weights as $s_{\mathrm{DG}}=\sum_j w_j U_j$, where $U_j$ are the nodal values on the element.

One-sided slopes are computed from neighboring cell averages,
\begin{equation}
d^- = \frac{\bar U - \bar U^{-}}{h},
\qquad
d^+ = \frac{\bar U^{+} - \bar U}{h},
\end{equation}
with the same construction applied in each coordinate direction in multiple dimensions.

The candidate and one-sided slopes are projected into characteristic space using the local eigenvectors of the Euler flux Jacobian,
\begin{equation}
A = L\,s_{\mathrm{DG}},
\qquad
B = L\,d^-,
\qquad
C = L\,d^+,
\end{equation}
where $L$ denotes the matrix of left eigenvectors evaluated at $\bar U$.

Each characteristic component is then limited using a TVB-modified minmod operator~\citep{cockburn1989tvb},
\begin{equation}
W_k =
\begin{cases}
A_k, & |A_k| \le M_{\mathrm{tvb}}\,h, \\[4pt]
\operatorname{sign}(A_k)
\min\!\left(|A_k|,\theta |B_k|,\theta |C_k|\right),
& \operatorname{sign}(A_k)=\operatorname{sign}(B_k)=\operatorname{sign}(C_k)\neq 0, \\[4pt]
0, & \text{otherwise},
\end{cases}
\end{equation}
where $M_{\mathrm{tvb}}$ is the empirically selected TVB parameter and $\theta>0$ is the minmod parameter.
The quantities $A_k$, $B_k$, and $C_k$ represent characteristic slopes, so the TVB criterion is applied directly to their magnitudes.

The limited characteristic slopes are mapped back to conserved-variable space as $s^{\mathrm{lim}}=R\,W$, where $R$ denotes the matrix of right eigenvectors.

The limiter is applied only to elements flagged by a KXRCF density-jump indicator~\citep{krivodonova2004shock}, while unflagged elements remain unchanged.
In multiple dimensions, the procedure is applied dimension by dimension using directional slopes and characteristic decompositions.

\section{Numerical experiments}
The numerical experiments in this section evaluate the IGR--DG method on a set of benchmark problems that reflect numerical challenges observed in shock-dominated compressible flows. 
The test cases span shock formation, propagation and interaction, contact discontinuities, and interactions with boundaries and flow structures.

The configurations are evaluated using two stabilization approaches.
The proposed IGR--DG method augments the Euler equations with the entropic pressure $\Sigma$, introducing stabilization at the PDE level.
As a baseline, the DG+limiter method is used as simple and well-understood reference to compare against the IGR--DG method.
The unstabilized DG discretizations are not included, since they typically develop nonphysical oscillations or singular behavior and do not provide a meaningful basis for comparison~\citep{gottlieb1997gibbs}.

The experiments assess stability as compressive intensity increases, the localization of the entropic pressure to regions of strong compression, and the preservation of fine scale flow features in the presence of shocks.
These properties are further examined in multidimensional problems where symmetry preservation and the absence of grid aligned artifacts provide additional measures of solution quality.
To quantify accuracy, mesh refinement studies are performed using the $L^{1}$ error to determine the convergence behavior of each method.
These results are complemented by qualitative comparisons that examine numerical diffusion near contact surfaces and the clipping of local extrema in smooth regions, both of which are commonly associated with limiter based methods~\citep{qiu2005runge}.
Because the IGR formulation requires smooth initial data, all test cases employ smooth initial conditions and classical Riemann type problems are represented using narrow smooth profiles~\citep{bhat2009regularization}.
This construction maintains consistency with the governing formulation while enabling direct comparison with conventional shock capturing methods.
The test cases cover flow regimes ranging from mildly varying acoustic fields to sharply localized transitions that rapidly develop into steep compressive structures.

\subsection{One-dimensional numerical experiments}
\label{sec:1d_experiments}

The one-dimensional numerical experiments examine the behavior of the IGR--DG method for a set of benchmark problems involving shock dominated flow.
These tests isolate the response of the regularization to compressive waves without the additional effects introduced by multidimensional flow.

Unless otherwise noted, all simulations are performed on $\Omega=[0,1]$ using piecewise-polynomial DG discretizations on uniform Cartesian grids with $m$ elements and explicit three-stage, third-order strong-stability-preserving Runge--Kutta SSPRK(3,3) time integration~\citep{gottlieb2001strong}.
Assuming a maximum velocity of $2$, a conservative Courant factor of $\mathrm{CFL}=0.4/(2p+1)$ is used~\citep{cockburn2001runge}, with time step $\Delta t=\mathrm{CFL}\,h/2$ unless otherwise specified.
The ratio of specific heats is fixed to $\gamma=1.4$, and an ideal compressible equation of state is used unless otherwise specified.
A local Lax--Friedrichs flux is used for the Euler equations, and the TVB parameter and minmod parameter are set to $M_{\mathrm{tvb}}=20$ and $\theta=1.5$, when the limiter is applied.

The regularization strength coefficient $\alpha_0$ is selected on a per-problem basis to match the level of compressive intensity.
Values are chosen to ensure stable evolution while avoiding unnecessary dissipation, and are varied in selected tests to examine their dependence on shock strength.
This parameter controls the strength and spatial scale of the regularization.

Discontinuous Riemann problems are initialized using smooth solution profiles, as required by the IGR formulation.
Left and right states are connected using a smooth hyperbolic tangent profile of fixed physical width, independent of grid resolution.
Specifically, a smooth step function
\begin{equation}
    s(x) = \tfrac12 \left( 1 - \tanh\!\left((x - x_0)/\delta\right) \right)
    \label{eq:smooth_step}
\end{equation}
is used to interpolate between left and right states across a transition region of width $\delta$ centered at $x_0$.
The initial condition is then defined by
\begin{equation}
    U(x,0) = s(x)\,U_L + \left(1 - s(x)\right)\,U_R,
    \label{eq:smooth_interpolate}
\end{equation}
where $U_L$ and $U_R$ denote the prescribed left and right states.
Holding the transition width $\delta$ fixed across mesh refinements avoids introducing grid-dependent initial conditions and ensures that observed convergence behavior reflects properties of the numerical scheme rather than changes in the prescribed initial structure.

The smooth transitions modify the classical Riemann problems, so the corresponding self-similar analytical solutions are no longer exact solutions of the prescribed initial value problems.
Unless otherwise specified, reference solutions are therefore obtained from fine-grid DG+limiter simulations.

\subsubsection{Shock formation and propagation}

    Baseline tests verify the correctness of the IGR--DG formulation using canonical one-dimensional shock dominated flows.
    These benchmarks distinguish between dynamically generated shock formation and sustained shock propagation.
    The nonlinear periodic simple wave problem represents the formation of shocks from smooth initial data, while the smoothed Sod shock tube represents the propagation of an established shock.

    \begin{example}[Nonlinear periodic simple wave]
        \label{ex:simple_wave}
        This example considers shock formation driven by nonlinear compression.
        The initial condition is a periodic simple wave with constant Riemann invariant, given by
        \begin{equation}
            \rho(x,0)= 1 + 0.2 \sin\left(\frac{2\pi x}{L}\right), \qquad
            u(x,0) = \frac{2}{\gamma-1}\left(c(\rho(x,0)) - c(1)\right),
        \end{equation}
        with pressure determined by the isentropic equation of state $P(\rho)=\rho^{\gamma}$, and $c(\rho)$ denotes the sound speed associated with the isentropic equation of state.
        Nonlinear steepening leads to shock formation at approximately $t=0.42$.
    \end{example}

    \begin{example}[Smoothed Sod shock tube]
        \label{ex:sod}
        This example examines shock propagation in a canonical one-dimensional configuration~\citep{sod1978survey}.
        The left and right primitive states are
        \begin{equation}
            (\rho_L,u_L,P_L) = (1.0,\,0.0,\,1.0), \qquad
            (\rho_R,u_R,P_R) = (0.125,\,0.0,\,0.1),
        \end{equation}
        with the discontinuity located at $x_0=L/2$.
        The transition between states is smoothed using the construction defined in \cref{eq:smooth_step,eq:smooth_interpolate} with fixed physical width $\delta = 3L/128$.
        Dirichlet boundary conditions are imposed at both domain boundaries.
        The solution consists of a left-moving rarefaction, a right-moving shock, and a contact discontinuity.
    \end{example}
    
\begin{figure}[pos=htbp]
        \centering
        % \tikzsetnextfilename{simple_wave_1D_snapshots_p3}
        % \input{figures/simple_wave_p3/figure_simple_wave_1D_snapshots_p3}
        \includegraphics{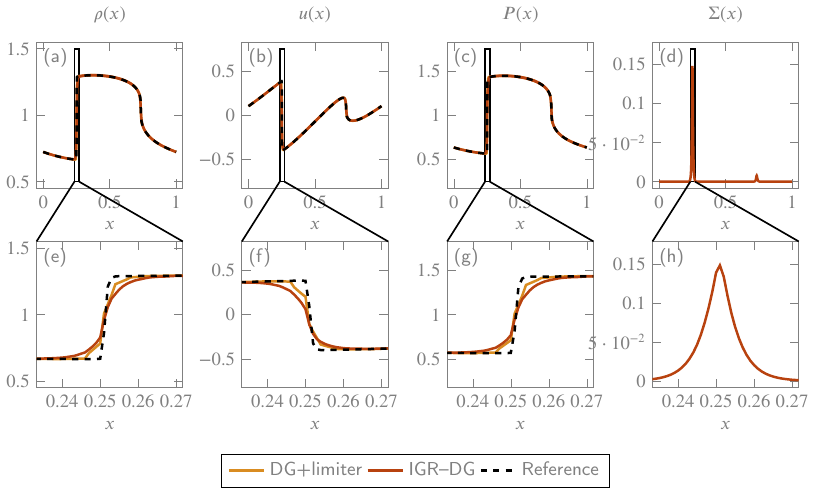}
        \caption{
            One-dimensional nonlinear simple-wave (Example~\ref{ex:simple_wave}) solution at polynomial order $p=3$.
            Numerical solutions are computed with $m=128$ cells and compared against a reference solution using $m_{\mathrm{ref}}=1024$ cells.
            Snapshots at $t=0.65$ are shown for (a,e) density $(\rho)$, (b,f) velocity $(u)$, (c,g) pressure $(P)$, and (d,h) entropic pressure $(\Sigma)$.
            Panels (a--d) show the full domain over $x\in[0,1]$, and panels (e--h) show the corresponding zoomed views over $x\in[0.235,0.27]$.
            Results are shown for DG+limiter and IGR--DG using $\alpha_{0}=0.3$.
        }
        \label{fig:simple_wave_snapshots_p3}
    \end{figure}

    In the nonlinear simple-wave problem, the shock forms dynamically from an initially smooth compressive wave.
    Snapshots at $t=0.65$ in \cref{fig:simple_wave_snapshots_p3} show that the entropic pressure $\Sigma$ is concentrated at the shock location and remains negligible elsewhere.
    Both the IGR--DG and DG+limiter methods capture the resulting shock without visible spurious oscillations.

\begin{figure}[pos=htbp]
        \centering
        % \tikzsetnextfilename{simple_wave_1D_error_multi_t_L1}
        % \input{figures/simple_wave_p3/figure_simple_wave_1D_error_multi_t_L1}
        \includegraphics{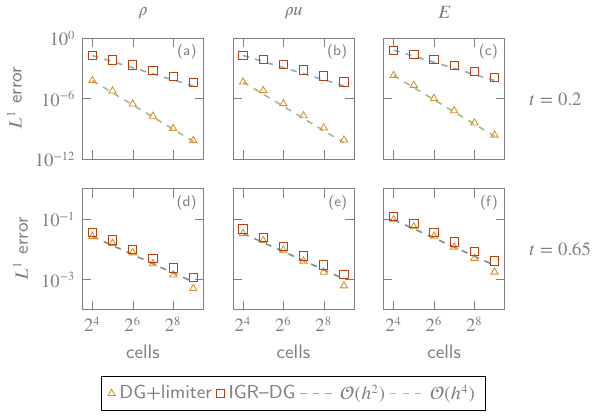}
        \caption{
            $L^{1}$ error convergence under grid refinement for the one-dimensional nonlinear simple-wave (Example~\ref{ex:simple_wave}) at polynomial order $p=3$.
            Errors are measured for (a,d) density $(\rho)$, (b,e) momentum $(\rho u)$, and (c,f) energy $(E)$.
            Panels (a--c) correspond to $t=0.2$, prior to shock formation, and panels (d--f) correspond to $t=0.65$, after shock formation.
            Results are shown for DG+limiter and IGR--DG using $\alpha_{0}=0.3$.
            Reference $\mathcal{O}(h^{2})$ and $\mathcal{O}(h^{4})$ scalings are shown at $t=0.2$, and reference $\mathcal{O}(h)$ scaling is shown at $t=0.65$.
        }
        \label{fig:simple_wave_errors_multi_t_L1}
    \end{figure}

    The error under grid refinement for the nonlinear simple-wave problem is reported in \cref{fig:simple_wave_errors_multi_t_L1}.
    Prior to shock formation, the DG+limiter method exhibits the expected fourth-order convergence for $p=3$, whereas the elliptic regularization reduces the observed convergence rate of the IGR--DG method to approximately second order.
    After shock formation, both methods exhibit first-order convergence.

\begin{figure}
        \centering
        % \tikzsetnextfilename{sod_1D_snapshots_p3}
        % \input{figures/sod_p3/figure_sod_1D_snapshots_p3}
        \includegraphics{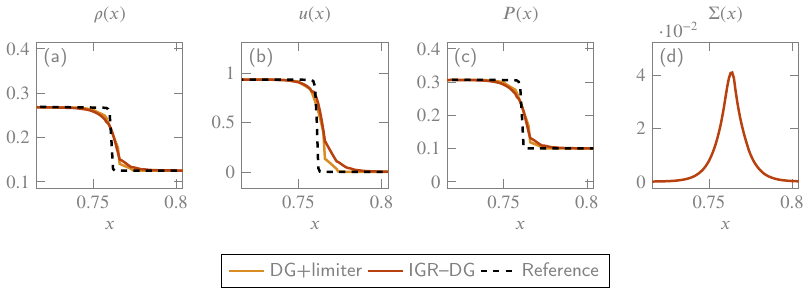}
        \caption{
            One-dimensional smoothed Sod shock tube (Example~\ref{ex:sod}) solution at polynomial order $p=3$.
            Numerical solutions are computed with $m=128$ cells and compared against a reference solution using $m_{\mathrm{ref}}=1024$ cells.
            Snapshots at $t=0.15$ are shown for (a) density $(\rho)$, (b) velocity $(u)$, (c) pressure $(P)$, and (d) entropic pressure $(\Sigma)$ over the region $x\in[0.72,0.8]$.
            Results are shown for DG+limiter and IGR--DG using $\alpha_{0}=0.7$.
        }
        \label{fig:sod_snapshots_p3}
    \end{figure}

\begin{figure}[pos=htbp]
        \centering
        % \tikzsetnextfilename{sod_1D_error_multi_t_L1}
        % \input{figures/sod_p3/figure_sod_1D_error_multi_t_L1}
        \includegraphics{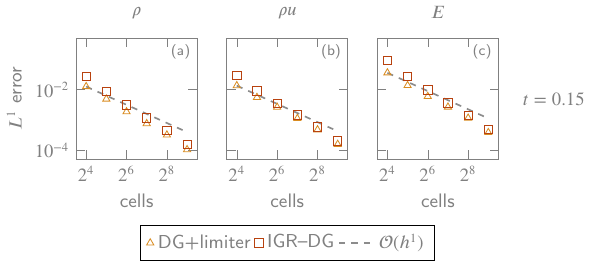}
        \caption{
            $L^{1}$ error convergence under grid refinement for the one-dimensional smoothed Sod shock tube (Example~\ref{ex:sod}) at polynomial order $p=3$.
            Errors are measured for (a) density $(\rho)$, (b) momentum $(\rho u)$, and (c) energy $(E)$ at $t=0.15$.
            Results are shown for DG+limiter and IGR--DG using $\alpha_{0}=0.5$.
        }
        \label{fig:sod_errors_multi_t_L1}
    \end{figure}

    In the smoothed Sod shock tube, the shock is present from the outset and propagates through the domain.
    Snapshots at $t=0.15$ in \cref{fig:sod_snapshots_p3} show accurate shock tracking and stable evolution for both methods.
    The error under grid refinement shown in \cref{fig:sod_errors_multi_t_L1} exhibits first-order convergence for both methods.

\subsubsection{Preservation of fine-scale features}

    These tests consider one-dimensional problems in which shocks interact with fine-scale features, emphasizing regimes where the strengths of the IGR--DG method are most pronounced.
    The primary challenge is to resolve both shocks and physically meaningful fine-scale features since limiter based methods often introduce dissipation that degrades high frequency features~\citep{dumbser2014posteriori}.
    The benchmarks presented here examine two settings, the modified Shu--Osher problem considers a propagating shock interacting with pre-existing oscillatory structure, while the high-frequency perturbation problem examines the evolution of fine-scale features during dynamically generated shock formation.

    \begin{example}[Modified Shu--Osher problem]
        \label{ex:shu_osher}
        This example based on the Shu--Osher problem~\citep{shu1989efficient}, examines the interaction between a propagating shock and pre-existing fine-scale features.
        The left and right primitive states are
        \begin{equation}
            (\rho_L,u_L,P_L) = (3.857143,\,2.629369,\,10.3333), \quad (\rho_R(x),u_R,P_R) = (\rho_b+\varepsilon_{\rho}\sin(kx),\,0.0,\,1.0),
        \end{equation}
        where $\rho_b = 1.0$, $\varepsilon_{\rho} = 0.2$, and $k = 50$.
        The discontinuity is located at $x_0=L/10$ and a sinusoidal density perturbation is superposed on the right state.
        The transition between states is smoothed using the construction defined in \cref{eq:smooth_step,eq:smooth_interpolate} with fixed physical width $\delta = 3L/128$.
        Dirichlet boundary conditions are imposed at both domain boundaries.
    \end{example}

    \begin{example}[High-frequency perturbation problem]
        \label{ex:high_freq_perturbation}
        This example complements the Shu--Osher problem by examining the retention of fine-scale features when a shock forms dynamically rather than being present initially.
        The configuration consists of a nonlinear compressive background wave with a superimposed high-frequency acoustic perturbation.
        The initial condition combines a low-frequency simple wave with a high-frequency pressure perturbation,
        \begin{equation}
            \begin{aligned}
                \rho(x,0) &= \rho_b + \varepsilon_{\rho} \sin\!\left(\frac{2\pi k_{1} x}{L}\right), \\
                u(x,0) &= \frac{2}{\gamma-1}\left(c(\rho(x,0)) - c(\rho_b)\right), \\
                P(x,0) &= P_{\mathrm{LF}}(x) + \varepsilon_{P} \sin\!\left(\frac{2\pi k_{2} x}{L}\right),
            \end{aligned}
        \end{equation}
        where \(P_{\mathrm{LF}}(x) = P(\rho(x,0))\) is the pressure associated with the low-frequency simple-wave state and \(c\) is the sound speed.
        The parameter values used in this work are \(\rho_b = 1.0\), \(\varepsilon_{\rho} = 0.25\), \(k_1 = 1\), \(\varepsilon_{P} = 0.05\), and \(k_2 = 16\).
        Periodic boundary conditions are imposed at both domain boundaries.
    \end{example}

\begin{figure}[pos=htbp]
        \centering
        % \tikzsetnextfilename{shu_osher_1D_snapshots_p3}
        % \input{figures/shu_osher_p3/figure_shu_osher_1D_snapshots_p3}
        \includegraphics{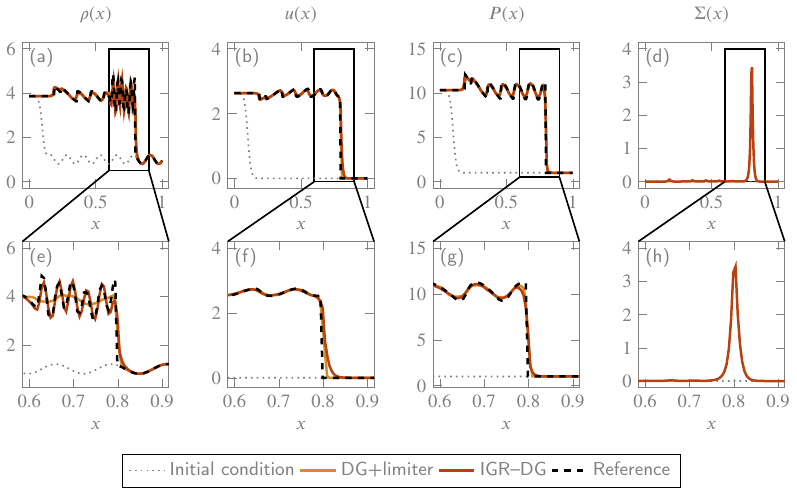}
        \caption{
            One-dimensional modified Shu--Osher problem (Example~\ref{ex:shu_osher}) at polynomial order $p=3$.
            Numerical solutions use $m=256$ cells and are compared with a reference solution using $m_{\mathrm{ref}}=1024$ cells.
            Snapshots at $t=0.2$ are shown for (a,e) density $(\rho)$, (b,f) velocity $(u)$, (c,g) pressure $(P)$, and (d,h) entropic pressure $(\Sigma)$.
            Panels (a--d) show the full domain $x\in[0,1]$, and panels (e--h) show the corresponding zoomed views over $x\in[0.6,0.9]$.
            Dotted curves indicate the initial condition, solid curves show DG+limiter and IGR--DG solutions using $\alpha_{0}=0.5$, and dashed curves indicate the reference solution.
        }
        \label{fig:shu_osher_snapshots_p3}
    \end{figure}

\begin{figure}[pos=htbp]
        \centering
        % \tikzsetnextfilename{shu_osher_1D_error_multi_t_L1}
        % \input{figures/shu_osher_p3/figure_shu_osher_1D_error_multi_t_L1}
        \includegraphics{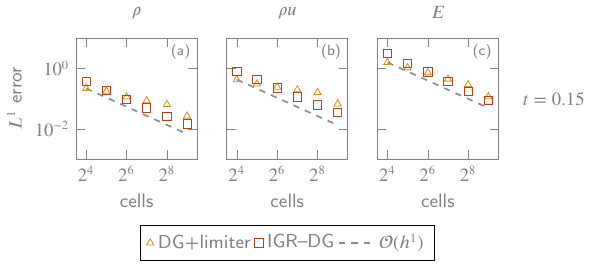}
        \caption{
            $L^{1}$ error convergence under grid refinement for the one-dimensional modified Shu--Osher problem (Example~\ref{ex:shu_osher}) at polynomial order $p=3$.
            Errors are measured for (a) density $(\rho)$, (b) momentum $(\rho u)$, and (c) energy $(E)$ at $t=0.15$.
            Results are shown for DG+limiter and IGR--DG using $\alpha_{0}=0.5$.
        }
        \label{fig:shu_osher_errors_multi_t_L1}
    \end{figure}

    In the modified Shu--Osher problem, a shock propagates through a pre-existing oscillatory density field.
    Snapshots at $t=0.2$ in \cref{fig:shu_osher_snapshots_p3} show that the IGR--DG method preserves the amplitude and phase of the post-shock oscillations more effectively than the DG+limiter method.
    Both methods capture the shock, but the DG+limiter method exhibits visible damping of downstream fine-scale features.
    The error under grid refinement at $t=0.15$, shown in \cref{fig:shu_osher_errors_multi_t_L1}, exhibits first-order convergence for both methods.
    The IGR--DG method yields lower error levels, particularly in density and total energy.

\begin{figure}[pos=htbp]
        \centering
        % \tikzsetnextfilename{high_freq_perturbation_1D_snapshots_p3}
        % \input{figures/high_freq_perturbation_p3/figure_high_freq_perturbation_1D_snapshots_p3}
        \includegraphics{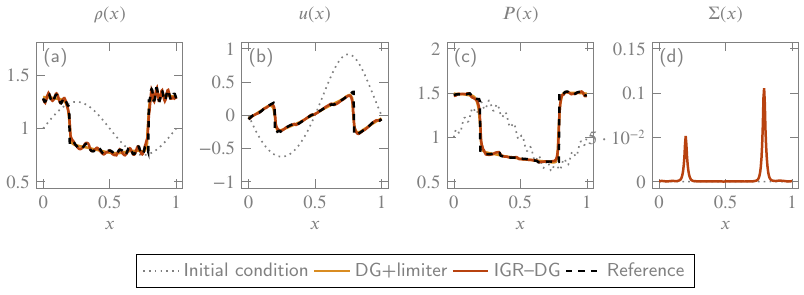}
        \caption{
            One-dimensional high-frequency perturbation problem (Example~\ref{ex:high_freq_perturbation}) at polynomial order $p=3$.
            Numerical solutions use $m=64$ cells and are compared with a reference solution using $m_{\mathrm{ref}}=1024$ cells.
            Snapshots at $t=1.0$ are shown for (a) density $(\rho)$, (b) velocity $(u)$, (c) pressure $(P)$, and (d) entropic pressure $(\Sigma)$ over the full domain $x\in[0,1]$.
            Dotted curves indicate the initial condition, solid curves show results for DG+limiter and IGR--DG using $\alpha_{0}=0.5$, and dashed curves indicate the reference solution.
        }
        \label{fig:hf_snapshots_p3}
    \end{figure}

\begin{figure}[pos=htbp]
        \centering
        % \tikzsetnextfilename{high_freq_perturbation_1D_error_multi_t_L1}
        % \input{figures/high_freq_perturbation_p3/figure_high_freq_perturbation_1D_error_multi_t_L1}
        \includegraphics{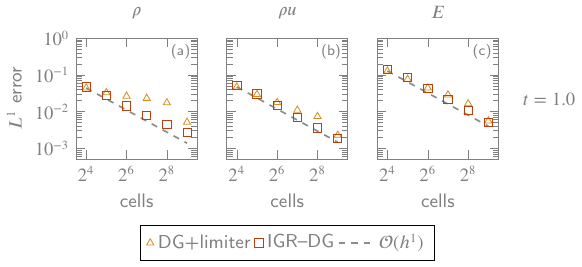}
        \caption{
            $L^{1}$ error convergence under grid refinement for the one-dimensional high-frequency perturbation problem (Example~\ref{ex:high_freq_perturbation}) at polynomial order $p=3$.
            Errors are measured for (a) density $(\rho)$, (b) momentum $(\rho u)$, and (c) energy $(E)$ at $t=1.0$.
            Results are shown for DG+limiter and IGR--DG using $\alpha_{0}=0.7$.
        }
        \label{fig:hf_errors_multi_t_L1}
    \end{figure}

    In the high-frequency perturbation problem, oscillatory structure evolves under increasing compression as a shock forms dynamically.
    Snapshots at $t=1.0$ in \cref{fig:hf_snapshots_p3} show that the oscillations persist through the forming shock, whereas the DG+limiter method progressively damps the high-frequency features.
    The error under grid refinement in \cref{fig:hf_errors_multi_t_L1} reflects this behavior, with the IGR--DG method maintaining lower error levels throughout the evolution.

    Across both problems, the IGR--DG method preserves fine-scale features more effectively while remaining stable in the presence of shocks.
    The regularization remains concentrated in strongly compressive regions, limiting damping away from shocks whether they are present initially or form dynamically during the evolution.

    \subsubsection{Benefit of higher-order approximation}

\begin{figure}[pos=htbp]
        \centering
        % \tikzsetnextfilename{shu_osher_fixed_dof256_1D_error_t_0-2_dof256}
        % \input{figures/shu_osher_fixed_dof256/figure_shu_osher_fixed_1D_error_t_0-2_dof256}
        \includegraphics{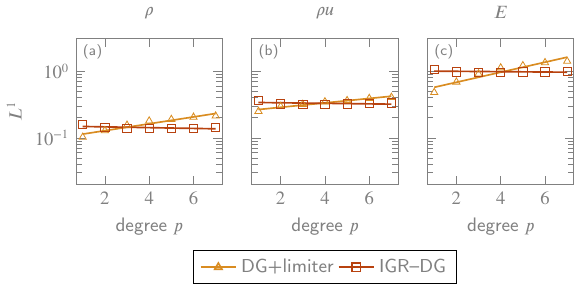}
        \caption{
            $L^{1}$ errors at fixed total degrees of freedom for the one-dimensional modified Shu--Osher problem (Example~\ref{ex:shu_osher}) under increasing polynomial order.
            Results are shown at $t=0.2$ with fixed $\mathrm{DOF}=256$.
            Errors are measured in density $(\rho)$, momentum $(\rho u)$, and energy $(E)$.
            Panels (a)--(c) show $\rho$, $\rho u$, and $E$, respectively.
            For each polynomial order $p$, the number of elements is chosen as $m \approx \mathrm{DOF}/(p+1)$.
            Results are shown for DG+limiter and IGR--DG with $\alpha_{0}=4.0$.
        }
        \label{fig:shu_osher_fixed_dof256_errors_t_0-2_L1}
    \end{figure}

    This study examines how solution accuracy varies with polynomial order when the number of degrees of freedom is held fixed, isolating the effect of $p$ refinement from $h$ refinement.
    The modified Shu--Osher problem is used because it combines shock propagation with fine-scale  features that remain sensitive to polynomial order.
    The total number of degrees of freedom is fixed at $\mathrm{DOF}=256$, and the number of elements is varied with $p$ according to $m \approx \mathrm{DOF}/(p+1)$.
    Time integration uses the SSPRK(3,3) scheme with $\mathrm{CFL}=0.4/(2p+1)$ and quadrature order $n_q=2p+3$~\citep{mengaldo2015dealiasing}.
    The regularization parameter is scaled as $\alpha=\alpha_{0}h_{\mathrm{eff}}^2$, where $h_{\mathrm{eff}}=2h/(p+1)$ denotes the effective nodal spacing within each element.
    This scaling maintains a consistent regularization length scale across polynomial orders.
    \Cref{fig:shu_osher_fixed_dof256_errors_t_0-2_L1} reports $L^{1}$ errors in density, momentum, and total energy at $t=0.2$.

    In the fixed DOF setting, the error is dominated by the shock, which reduces the effective approximation order, so increasing $p$ does not reduce the error.
    As $p$ increases, the number of elements decreases, further limiting resolution of the shock and post-shock features.
    Although increasing the polynomial degree improves the approximation properties of DG without modifying the mesh and enables finer scales to be resolved when the solution remains sufficiently smooth~\citep{ainsworth2004dispersive}, these gains can be diminished in shock dominated flows when shock capturing mechanisms introduce dissipation that suppress higher order modes~\citep{trojak2021shock}.

    For the DG+limiter method, the error increases with $p$, consistent with traditional limiter activation near the shock and in the downstream oscillatory region.
    The TVB limiter suppresses higher order modes in regions containing fine scale features, leading to clipping of smooth extrema and preventing higher polynomial orders from improving resolution even in well resolved parts of the domain~\citep{krivodonova2007limiters}.
    This effect becomes more pronounced in the fixed DOF setting because increasing $p$ necessarily reduces the number of elements available to resolve the shock and post shock features.

    In contrast, IGR--DG exhibits nearly constant error across polynomial orders, indicating that the regularization stabilizes the solution without a corresponding increase in shock dominated error as $p$ increases.
    This behavior is consistent with the motivation behind localized continuous regularization approaches, which seek to preserve the underlying high-order approximation while providing shock stabilization~\citep{bai2022continuous}.

\begin{figure}[pos=htbp]
        \centering
        % \tikzsetnextfilename{shu_osher_fixed_dof256_1D_snapshots_dof256}
        % \input{figures/shu_osher_fixed_dof256/figure_shu_osher_1D_snapshots_dof256}
        \includegraphics{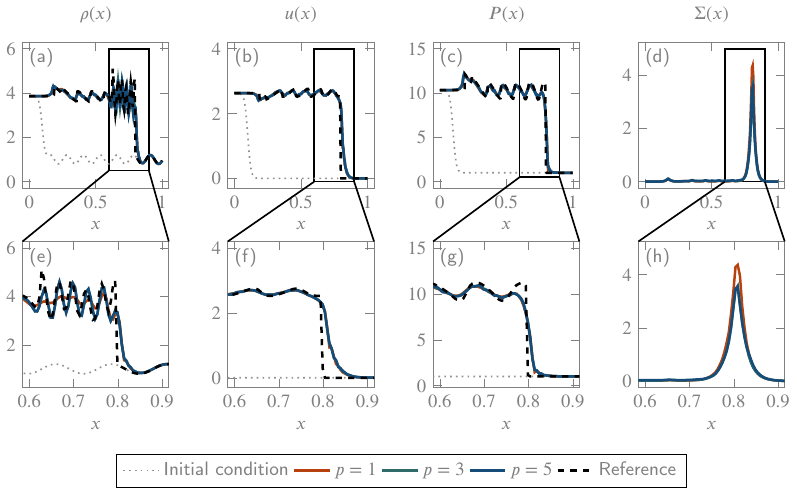}
        \caption{
            One-dimensional modified Shu--Osher problem (Example~\ref{ex:shu_osher}) solution at fixed $\mathrm{DOF}=256$.
            Numerical solutions are computed using IGR--DG with $\alpha_{0}=4.0$ for polynomial orders $p=1,3,5$, corresponding to $m=128,64,43$ elements, and compared against an IGR--DG reference solution using $m_{\mathrm{ref}}=1024$ elements, $p_{\mathrm{ref}}=3$, and $\alpha_{0,\mathrm{ref}}=1.0$.
            Snapshots are shown for (a,e) density $(\rho)$, (b,f) velocity $(u)$, (c,g) pressure $(P)$, and (d,h) entropic pressure $(\Sigma)$.
            Panels (a--d) show the full domain at $t=0.2$, and panels (e--h) show the corresponding zoomed views of the post-shock region.
            Dotted curves indicate the initial condition, and solid curves show results for different polynomial orders.
        }
        \label{fig:shu_osher_snapshots_dof256}
    \end{figure}

    \cref{fig:shu_osher_snapshots_dof256} shows solution snapshots at $t=0.2$ for the modified Shu--Osher problem for various polynomial orders.
    Differences are localized to the oscillatory region behind the shock, where the $p=1$ solution smooths the features while the $p=5$ solution retains higher-frequency content.
    Away from this region, the solutions are similar across polynomial orders, with differences becoming more apparent at low DOF where limited resolution and reduced ability of lower $p$ to represent curvature diminish fine-scale features.

\subsubsection{Stability under increasing Mach number}

    This subsection evaluates how the minimum regularization strength, $\alpha$, required for stable evolution varies with compressive intensity.
    To characterize this dependence, compressive intensity is selected by varying the peak Mach number while all other aspects of the setup remain fixed.
    As the initial Mach number increases, the solution transitions from smooth acoustic propagation to strongly nonlinear evolution with shock formation, requiring progressively stronger regularization to suppress the growth of numerical oscillations and prevent the simulation from failing.
    Stable evolution is recovered when $\alpha$ exceeds a threshold, and this minimum stabilizing value of $\alpha$ increases with shock strength.

    \begin{example}[Sinusoidal compressive wave]
        \label{ex:sinusoid_mach}
        This example considers the evolution of a compressive wave with fixed spatial structure and variable amplitude.
        The initial condition is parameterized by the target Mach number $M_\star$:        
        \begin{equation}
            \begin{aligned}
                \rho(x,0) &= 1.0, \\
                u(x,0) &= M_\star c_0 \sin\left(2 \pi k x\right), \\
                P(x,0) &= 1 + \delta_e \sin\left(2 \pi k x + \psi\right),
            \end{aligned}
        \end{equation}
        where $c_0 = c(1.0)$ is the sound speed associated with the base state.
        The parameter $M_\star$ controls the peak Mach number of the initial condition through the velocity amplitude. 
        By varying $M_\star$, this configuration isolates the effect of increasing compressive intensity while maintaining a fixed spatial structure.
        For small values of $M_\star$, the solution remains smooth, while larger values lead to nonlinear steepening and the formation of shock-like structures.
        The parameter values used in this work are $\delta_e = 0.20$, $\psi = \frac{\pi}{2}$, and $k = 1$.
        Periodic boundary conditions are imposed at both domain boundaries.
    \end{example}

\begin{figure}[pos=htbp]
        \centering
        % \tikzsetnextfilename{mach_alpha_sinusoid_1D_failure_point}
        % \input{figures/mach_sinusoid_p1/figure_mach_alpha_sinusoid_1D_failure_point}
        \includegraphics{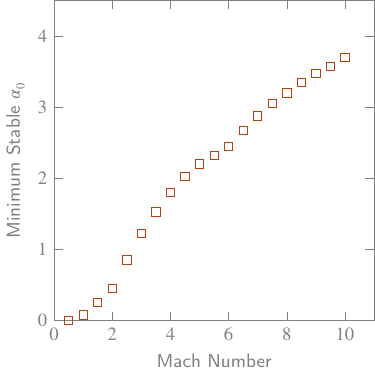}
        \caption{
            Minimum stabilizing value of $\alpha_{0}$ for the one-dimensional sinusoidal compressive wave (Example~\ref{ex:sinusoid_mach}) with polynomial order $p=1$.
            Markers show the smallest $\alpha_{0}$ yielding stable evolution as a function of the target Mach number $M_\star$.
            Stability is assessed over $M_\star\in[0.5,10.0]$, with $\alpha_{0}$ sampled in increments of $0.025$.
        }
        \label{fig:stability_boundary_alpha_mach}
    \end{figure}

    The dependence of the minimum stabilizing value of $\alpha$ on $M_\star$ for linear basis functions ($p=1$) is shown in \cref{fig:stability_boundary_alpha_mach}.
    As the target $M_\star$ increases, the minimum stabilizing value of $\alpha$ also increases, reflecting the stronger regularization required as compressive gradients become steeper.
    Across the range of $M_\star$ considered, this threshold lies within $\alpha_{0}\in[0.5,\,4.0]$ for the linear basis functions ($p=1)$, although larger values may be used to suppress residual oscillations.
    For higher polynomial orders, the same trend is expected after accounting for the dependence of the regularization on $p$.
    Generally, $\alpha$ should scale with the magnitude of compressive gradients, providing a physically interpretable basis for parameter selection rather than ad hoc tuning.

    \begin{figure}[pos=htbp]
        \centering
        % \tikzsetnextfilename{mach_alpha_sinusoid_1D_snapshot_Mach4}
        % \input{figures/mach_sinusoid_p1/figure_mach_alpha_sinusoid_1D_snapshot}
        \includegraphics{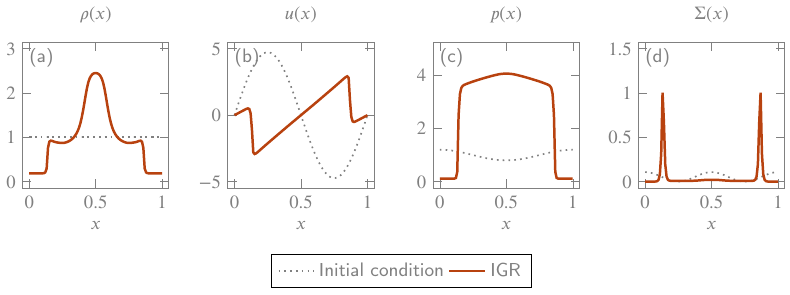}
        \caption{
            One-dimensional sinusoidal compressive wave (Example~\ref{ex:sinusoid_mach}) at target Mach number $M_\star=4.0$ with polynomial order $p=1$ and $m=256$ elements.
            Snapshots at $t=0.15$ are shown for (a) density $(\rho)$, (b) velocity $(u)$, (c) pressure $(P)$, and (d) entropic pressure $(\Sigma)$.
            Dotted curves indicate the initial condition, and solid curves show the IGR--DG solution using $\alpha_{0}=4.0$.
        }
        \label{fig:sinusoid_mach4_snapshots}
    \end{figure}

\subsubsection{Shock interactions and boundary effects}

    This subsection considers one-dimensional problems in which shocks interact with boundaries or with other shocks.
    These settings place increased demands on robustness because compressive structures can reflect and overlap within confined domains, producing multiple regions of strong compression.
    The benchmark problems distinguish between interactions with physical boundaries and interactions between shocks, allowing these two mechanisms to be examined separately.
    Shock reflection at a rigid wall isolates the interaction between a compressive wave and a boundary, while the double blast wave problem examines the interaction of multiple shocks within a bounded domain.

    \begin{example}[Shock reflection at a rigid wall]
        \label{ex:shock_reflection}
        This example examines the interaction of a shock with a reflective boundary.
        The initial condition is defined by left and right primitive states
        \begin{equation}
            (\rho_L,u_L,P_L) = (3.0,\,2.0,\,3.0), \qquad
            (\rho_R,u_R,P_R) = (1.0,\,0.0,\,1.0),
        \end{equation}
        with the discontinuity located at $x_0 = L/3$.
        The initial condition is constructed using the smoothing procedure defined in \cref{eq:smooth_step,eq:smooth_interpolate} with fixed physical width $\delta = 3L/128$.
        Reflective boundary conditions are imposed at the right boundary, with Dirichlet conditions at the left boundary.
        This shock reflection configuration produces a right-moving shock that reflects at the boundary, generating complex wave interactions.
        It isolates the ability of the method to handle strong compression and boundary-induced reflections without introducing nonphysical oscillations.
    \end{example}

    \begin{example}[Double blast wave]
        \label{ex:double_blast_wave}
        This example examines shock--shock interaction generated by two pressure perturbations in a bounded domain.
        The configuration produces inward-moving shocks and outward-moving rarefaction waves that interact nonlinearly near the center of the domain.
        The initial condition is defined by three constant primitive states
        \begin{equation}
            (\rho_L,u_L,P_L) = (1.0,\,0.0,\,3.0), \qquad
            (\rho_M,u_M,P_M) = (1.0,\,0.0,\,1.0), \qquad
            (\rho_R,u_R,P_R) = (1.0,\,0.0,\,3.0),
        \end{equation}
        separated by interfaces at $x_L=3L/10$ and $x_R=7L/10$.
        The transitions between the three states are smoothed using the construction defined in \cref{eq:smooth_step,eq:smooth_interpolate} with fixed physical width $\delta = 3L/128$.
        Reflective boundary conditions are imposed at both domain boundaries.
        The resulting blast waves propagate inward and collide near the center of the domain, producing interacting shocks, contact layers, and reflected waves.
    \end{example}

\begin{figure}[pos=htbp]
        \centering
        % \tikzsetnextfilename{shock_reflection_1D_snapshots_p3}
        % \input{figures/shock_reflection_p3/figure_shock_reflection_1D_snapshots_p3}
        \includegraphics{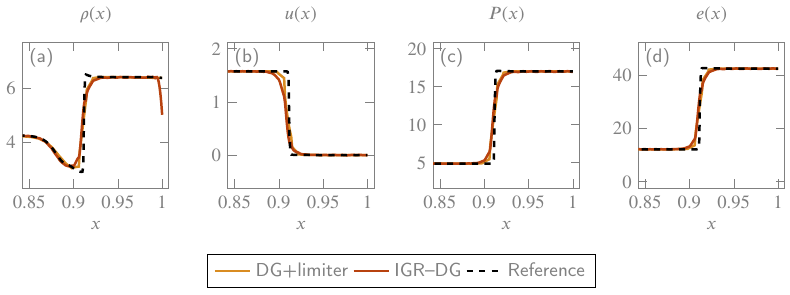}
        \caption{
            One-dimensional shock reflection problem (Example~\ref{ex:shock_reflection}) at polynomial order $p=3$ with $m=128$ elements.
            Snapshots at $t=0.35$ are shown for (a) density $(\rho)$, (b) velocity $(u)$, (c) pressure $(P)$, and (d) internal energy $(e)$ over the region near the reflected shock.
            Solid curves show DG+limiter and IGR--DG solutions with $\alpha_{0}=0.5$, and dashed curves indicate the reference solution computed with $m_{\mathrm{ref}}=1024$ elements.
        }
        \label{fig:shock_reflection_snapshots_p3}
    \end{figure}

    In the shock reflection problem, a single strong shock interacts with a rigid wall.
    Snapshots at $t=0.35$ in \cref{fig:shock_reflection_snapshots_p3} show that both the DG+limiter and IGR--DG methods capture the reflected shock without visible spurious oscillations.
    The IGR--DG solution exhibits a small mismatch between boundary-adjacent and interior values of density and internal energy after reflection, consistent with the wall heating effect commonly observed in numerical simulations of strong shocks~\citep{noh1987errors}.
    This effect remains confined near the wall and does not noticeably affect the interior solution.

\begin{figure}[pos=htbp]
        \centering
        % \tikzsetnextfilename{double_blast_wave_1D_snapshots_p3}
        % \input{figures/double_blast_wave_p3/figure_double_blast_wave_1D_snapshots_p3}
        \includegraphics{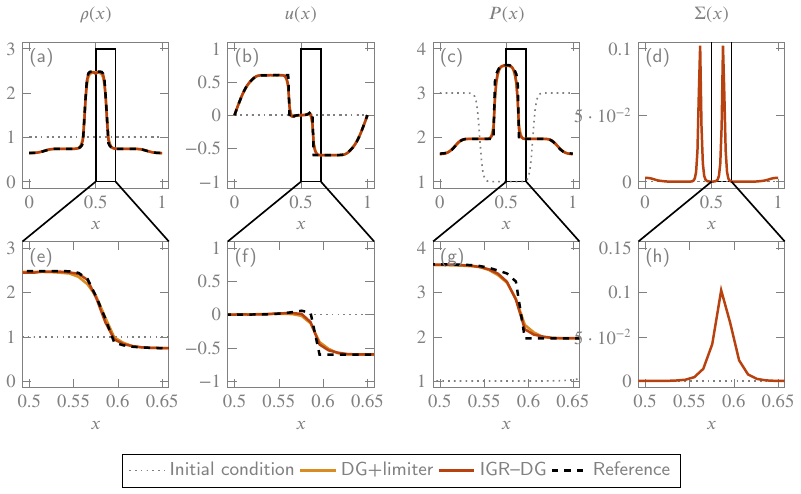}
        \caption{
            One-dimensional double blast wave problem (Example~\ref{ex:double_blast_wave}) at polynomial order $p=3$ with $m=64$ elements.
            Snapshots at $t=0.2$ are shown for (a,e) density $(\rho)$, (b,f) velocity $(u)$, (c,g) pressure $(P)$, and (d,h) entropic pressure $(\Sigma)$.
            Panels (a--d) show the full domain, and panels (e--h) show a zoomed region near the interacting shocks.
            Solid curves show DG+limiter and IGR--DG solutions with $\alpha_{0}=0.3$, dashed curves indicate the reference solution computed with $m_{\mathrm{ref}}=1024$ elements, and dotted curves show the initial condition.
        }
        \label{fig:double_blast_wave_snapshots_p3}
    \end{figure}

    In the double blast wave problem, interacting shocks create overlapping regions of strong compression within the domain.
    Snapshots at $t=0.2$ in \cref{fig:double_blast_wave_snapshots_p3} show that the entropic pressure $\Sigma$ forms localized peaks aligned with the shock fronts.
    The IGR--DG method produces sharper density profiles and more closely follows the reference solution in the vicinity of the shocks.

\subsection{Two-dimensional numerical experiments}

The two-dimensional numerical experiments examine the behavior of the IGR--DG method in flows involving multidirectional compression, shock interaction, and complex wave propagation.
Relative to the one-dimensional benchmarks, these configurations introduce geometric effects that require the regularization to remain localized while preserving symmetry and avoiding grid-aligned artifacts.
The test cases considered here include radially symmetric blast waves, shock interactions with interfaces and vortical structures, and compressible flows around curved boundaries.

\subsubsection{Julia implementation}

    The following experiments use a Julia-based DG solver obtained by extending the one-dimensional implementation of \Cref{sec:1d_experiments} to two spatial dimensions.
    This version is designed for simplicity and accessibility, allowing the method to be implemented, modified, and executed without reliance on an external CFD framework.
    The examples extend the one-dimensional studies of stability under compression and shock interaction to multidimensional flows.
    The source code for the Julia implementation used in this subsection is available in the \texttt{igr\_dg\_julia} repository~\citep{eyob2026igrdgjulia}.

\subsubsection*{Stability under strong compression}

    A smoothed 2D blast wave example is used to examine stability and symmetry preservation under strong multidimensional compression.
    The radial symmetry of the configuration makes it sensitive to directional bias and grid-aligned artifacts.

    \begin{example}[Smoothed 2D blast wave]
        \label{ex:blast_wave_2d}
        This example considers a two-dimensional blast wave generated by a localized deposition of internal energy.
        The problem is posed on $\Omega=[0,1]\times[0,1]$ with reflective boundary conditions.
        The background state is uniform, with $\rho=1.0$, $\rho\vct{u}=(0,0)^\top$, and $P=1.0$.

        To obtain a smooth initial condition, the blast is initialized through a radially symmetric pressure perturbation centered at $x_c=(\tfrac12,\tfrac12)$.
        Given
        \begin{equation}
            r(x)=\sqrt{(x_1-\tfrac12)^2+(x_2-\tfrac12)^2},
        \end{equation}
        the initial pressure is
        \begin{equation}
            P(x,0)=P+\epsilon\exp\!\left(-(r/\delta)^2\right),
        \end{equation}
        where
        \begin{equation}
            \delta=\frac{6}{32}\min(L_1,L_2), \text{ and }
            \qquad\qquad
            \epsilon=\frac{(\gamma-1)E_{\mathrm{blast}}}
            {\displaystyle\int_{\Omega}\exp\!\left(-(r/\delta)^2\right)\,dx},
        \end{equation}
        so that the total added internal energy equals $E_{\mathrm{blast}}=1.0$.

        The solution consists of a radially symmetric, outward-propagating shock separating a high-pressure interior from the ambient state.
    \end{example}

\begin{figure}[pos=htbp]
        \centering
        % \tikzsetnextfilename{blast_wave_2d_julia_snapshots}
        % \input{figures/blast_wave_2d/blast_wave_2d_julia_snapshots}
        \includegraphics{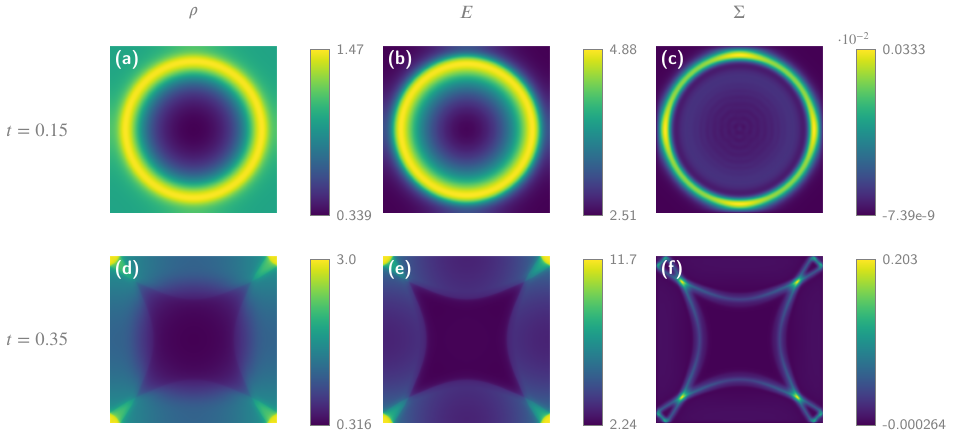}
        \caption{
            Two-dimensional blast wave problem (Example~\ref{ex:blast_wave_2d}) at polynomial order $p=3$ on a $128\times128$ mesh.
            Snapshots at $t=0.15$ (top row) and $t=0.35$ (bottom row) are shown for (a,d) density $(\rho)$, (b,e) total energy $(E)$, and (c,f) entropic pressure $(\Sigma)$.
            Results are shown for IGR--DG using $\alpha_{0}=5.0$.
        }
        \label{fig:blast_wave_2d_julia_snapshots}
    \end{figure}

    Solution snapshots for the 2D blast wave are shown in \cref{fig:blast_wave_2d_julia_snapshots} at $t=0.15$ and $t=0.35$.
    The initial energy deposition generates a radially expanding wave that steepens into an outward-propagating shock.
    At $t=0.15$, a narrow transition layer has formed at the shock front, separating a high-energy interior from the ambient state.
    By $t=0.35$, the shock has reached the domain boundaries, where it reflects while remaining well resolved.
    The entropic pressure $\Sigma$ is concentrated near the shock front and remains negligible elsewhere, exhibiting minimal response in both the interior and exterior regions.
    This demonstrates that the regularization activates selectively in regions of strong compression while remaining inactive elsewhere.
    Although the entropic pressure is predominantly positive, small localized negative values of $\Sigma$ are also observed, reflecting the contribution of non-compressive deformations to the right-hand side of the elliptic equation.

    The density and energy fields remain radially symmetric prior to boundary interaction, with no visible grid-aligned artifacts.
    Following reflection, the solution preserves the expected symmetry of the square domain, and the shock maintains a consistent thickness throughout the evolution.
    The solutions are free of oscillations and preserve the expected symmetry, demonstrating stable behavior under strong multidimensional compression.

\subsubsection*{Shock--structure interaction}

    A shock--bubble interaction example is used to examine the behavior of the regularization where the shock interacts with an interface.
    The focus is on the resolution of interface deformation and the localization of the regularization in the presence of coupled compressive and shear dynamics.

    \begin{example}[Shock--bubble interaction]
        \label{ex:shock_bubble_julia}
        A planar shock interacts with a circular bubble separating regions of differing density.
        The problem is posed on $\Omega=[0,2]\times[0,1]$ with Dirichlet boundary conditions on the left, transmissive conditions on the right, and reflective boundaries on the top and bottom.
        The background state has density $\rho_{\mathrm{bg}}=10.0$, velocity $u=(0,0)$, and pressure $P_{\mathrm{amb}}=1.0$.
        A circular bubble is centered at $x_c=(1,0.5)$ with radius $R_b=0.1$.
        The density inside the bubble is $\rho_{\mathrm{bub}}/\rho_{\mathrm{bg}}=\chi$ with $\chi=20.0$.
        The transition between the bubble and background is smoothed using the construction defined in \cref{eq:smooth_step,eq:smooth_interpolate} with smoothing width $\delta_b$.
        A planar shock is initialized at $x_{\mathrm{shock}}=1.0$.
        The transition between the pre- and post-shock states is smoothed using the construction defined in \cref{eq:smooth_step,eq:smooth_interpolate} with smoothing width $\delta_s$.
        The pre- and post-shock primitive states are
        \begin{equation}
            (\rho_{\mathrm{pre}},u_{x,\mathrm{pre}},u_{y,\mathrm{pre}},P_{\mathrm{pre}}) = (10.0,\,0.0,\,0.0,\,1.0), \qquad
            (\rho_{\mathrm{post}},u_{x,\mathrm{post}},u_{y,\mathrm{post}},P_{\mathrm{post}}) = \left(\frac{80}{3},\,0.4677,\,0.0,\,4.5\right),
        \end{equation}
        corresponding to a Mach number $M_s=2.0$ through the Rankine--Hugoniot relations.
        The initial condition is constructed by combining the smoothed shock profile with the spatially varying density field induced by the bubble, producing a smooth approximation of a shock interacting with a density interface.
    \end{example}

\begin{figure}[pos=htbp]
        \centering
        % \tikzsetnextfilename{shock_bubble_2d_julia_snapshots}
        % \input{figures/shock_bubble_2d/shock_bubble_2d_julia_snapshots}
        \includegraphics{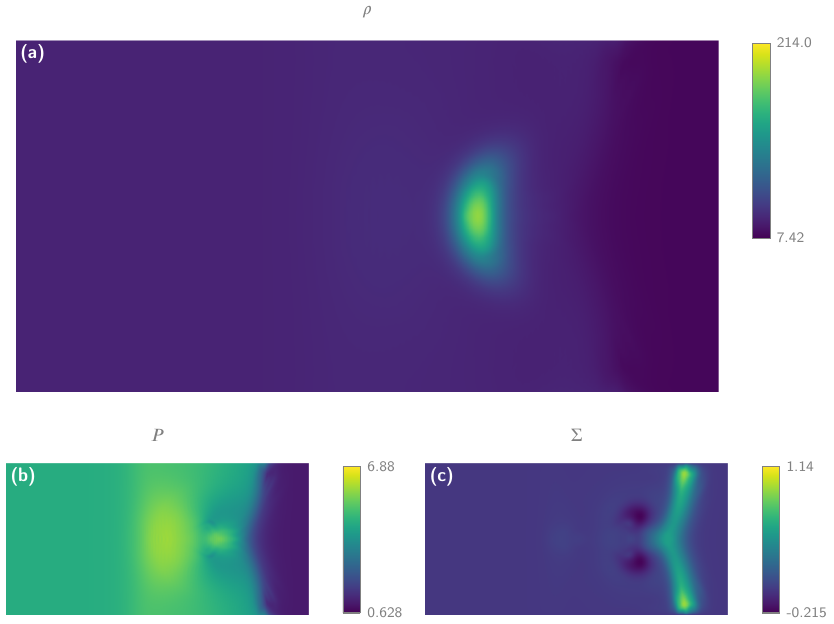}
        \caption{
            Two-dimensional shock--bubble interaction problem (Example~\ref{ex:shock_bubble_julia}) at polynomial order $p=3$ on a $128\times64$ mesh.
            The initial condition uses smoothing widths $\delta_b=3/64$ for the bubble interface and $\delta_s=3/32$ for the shock.
            Snapshots at $t=1.0$ are shown for (a) density $(\rho)$, (b) pressure $(P)$, and (c) entropic pressure $(\Sigma)$.
            Results are shown for IGR--DG using $\alpha_{0}=32.0$.
        }
        \label{fig:shock_bubble_2d_snapshots}
    \end{figure}

    Solution snapshots for the shock--bubble interaction at time $t=1.0$ are shown in \cref{fig:shock_bubble_2d_snapshots}.
    At this stage, the planar shock has propagated through the bubble and become curved, advancing more rapidly away from the centerline while deforming the density interface.
    The initially circular bubble is compressed and elongated in the streamwise direction while remaining symmetric about the horizontal centerline, matching the expected behavior of a shock interacting with a density interface.
    The pressure field captures the resulting wave interactions, while the main structure of the bubble is preserved throughout the deformation.

    The entropic pressure $\Sigma$ is localized near the shock and regions of strong compression.
    For this problem, a larger regularization parameter, $\alpha_{0}=32.0$, is used to accommodate the stronger compressive dynamics generated during the shock--bubble interaction.
    Away from these regions, the response remains small, indicating that the regularization activates selectively rather than acting as a uniform dissipative mechanism along the interface.
    No grid-aligned artifacts are observed, and the symmetry of the solution is preserved.
    These results show that the IGR--DG method captures the principal features of shock--interface interaction while maintaining localized, non-diffusive regularization and stable multidimensional behavior.

\subsubsection{MFEM implementation}

An MFEM implementation is presented here that provides an independent realization of the IGR--DG formulation within a scalable finite element framework.
In contrast to the structured Cartesian meshes used in the Julia implementation, this solver supports unstructured meshes, curved geometries, and more general boundary conditions.
The examples examine whether the behavior observed in the Julia-based benchmarks persists in settings with greater geometric and numerical complexity.
They focus on multidimensional shock interactions involving fine-scale features and curved boundaries that are not present in the structured Cartesian test cases.
The source code for the MFEM implementation used in this subsection is available in the pinned \texttt{dg-igr-mfem} reproduction repository~\citep{arias2026dgigrmfem}.

All MFEM results shown in this subsection use a discontinuous Galerkin finite element discretization of the Euler equations, with numerical fluxes across interelement faces and element-local inversion of the DG mass matrix~\citep{mfem-2021}.
The accompanying IGR problem is posed on an $L^2$ discontinuous space and discretized using a SIPG method.
Time integration uses the SSPRK(3,3) together with CFL-based adaptive timestepping.

\subsubsection*{Stability under strong compression}

This example revisits the 2D blast wave configuration introduced in \cref{ex:blast_wave_2d} to assess the consistency of the IGR--DG formulation across different discretization frameworks.
While the Julia implementation is based on a structured Cartesian mesh with tensor-product basis functions, the MFEM solver employs a general finite element infrastructure capable of handling unstructured meshes and more complex geometries.
The physical setup matches the Julia-based 2D blast wave configuration, with the blast initialized on $\Omega=[0,1]\times[0,1]$ from the same Gaussian energy deposition centered at $(\tfrac12,\tfrac12)$ and slip-wall boundaries imposed on all sides.

\begin{figure}[pos=htbp]
        \centering
        % \tikzsetnextfilename{blast_wave_2d_mfem_snapshots}
        % \input{figures/blast_wave_2d/blast_wave_2d_mfem_snapshots}
        \includegraphics{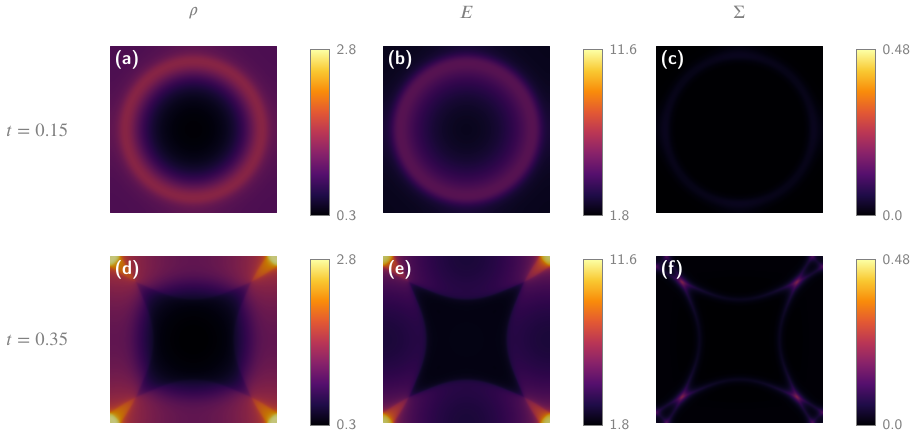}
        \caption{
            Two-dimensional blast wave problem (Example~\ref{ex:blast_wave_2d}) at polynomial order $p=3$.
            Snapshots are shown for (a,d) density $(\rho)$, (b,e) total energy $(E)$, and (c,f) entropic pressure $(\Sigma)$.
            Rows correspond to times $t=0.15$ and $t=0.35$ from top to bottom.
            Results are shown for IGR--DG using $\alpha_{0}=5.0$.
        }
        \label{fig:blast_wave_2d_snapshots}
    \end{figure}

The MFEM results shown in \cref{fig:blast_wave_2d_snapshots} use a $32\times 32$ base mesh, polynomial order $p=3$, two levels of parallel refinement, CFL$=0.25$, $\alpha_{0}=5.0$, and the $L^2$ SIPG formulation described above.
Despite these implementation differences, the MFEM results closely match those obtained with the Julia-based solver.
The entropic pressure activates in the vicinity of the blast wave, remains localized to the expanding shock front, and produces a smooth, radially symmetric solution profile.
No spurious oscillations or grid-aligned artifacts are observed.
This agreement indicates that the behavior of the IGR regularization is not tied to a specific discretization strategy and instead reflects a property of the underlying PDE-level formulation.

\subsubsection*{Shock interaction with fine-scale features}

The Inoue--Hattori shock--vortex interaction is used to examine the behavior of the regularization in the presence of both compressive and rotational dynamics.
The regularization must distinguish between these mechanisms, activating near the shock while remaining minimally dissipative within the vortex.
The shock is not aligned with the mesh and tests both robustness under grid non-alignment and preservation of fine-scale features through shock interaction.
This problem extends the one-dimensional Shu--Osher test to two dimensions, with vorticity replacing oscillatory density structure.

\begin{example}[Inoue--Hattori shock--vortex interaction]
    \label{ex:shock_vortex}
    The MFEM shock--vortex case is based on the large-domain Inoue--Hattori configuration introduced in~\citep{inoue1999sound} on $\Omega=[-30,10]\times[-20,20]$, discretized by a $32\times 32$ Cartesian quadrilateral base mesh.
    The initial condition combines a right-to-left normal shock of Mach number $M_s=2.0$ located at $x_s=-10$ with an isentropic vortex of strength $M_v=0.25$ centered at $(x_v,y_v)=(5,0)$.
    Relative to the original Inoue--Hattori configuration, the shock is smoothed over a few mesh cells and both the shock and vortex locations are shifted to fit the computational window used here.
    The shock profile is smoothed with a hyperbolic-tangent transition of width $\delta = 4\min\{h\}$, where $h$ denotes the mesh cell size in the two coordinate directions.
    A supersonic inflow state is imposed at the right boundary, a non-reflecting outlet state is imposed at the left boundary, and the top and bottom boundaries are treated as slip walls.
\end{example}

\begin{figure}[pos=htbp]
    \centering
    % \tikzsetnextfilename{shock_vortex_2d_snapshots_10-0}
    % \input{figures/shock_vortex_2d/shock_vortex_2d_snapshots_10-0}
    \includegraphics{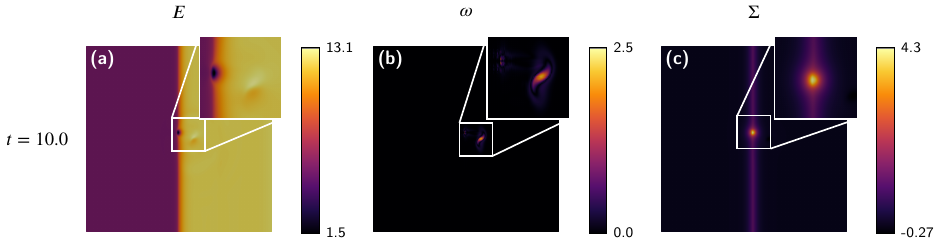}
    \caption{
        Two-dimensional Inoue--Hattori shock--vortex interaction problem (Example~\ref{ex:shock_vortex}) at polynomial order $p=3$.
        Snapshots at $t=10.0$ are shown for (a) total energy $(E)$, (b) vorticity $(\omega)$, and (c) entropic pressure $(\Sigma)$.
        Zoomed views of the shock--vortex interaction region are shown in the upper-right corner of each panel.
        Results are shown for IGR--DG using $\alpha_{0}=5.0$ and $\mathrm{CFL}=0.25$.
    }
    \label{fig:shock_vortex_2d_snapshots}
\end{figure}

As the solution evolves, the vortex convects from right to left toward the stationary shock and interacts with it.
Representative snapshots of the total energy, vorticity, and entropic pressure at $t=10.0$ are shown in \cref{fig:shock_vortex_2d_snapshots}.
The vortex is compressed and stretched as it passes through the shock, producing filamentary structures in the vorticity field due to differential advection and compression.
The shock remains sharply defined, with $\Sigma$ activating along the compressive interface.
A localized energy deficit appears along the shock at the vortex interaction point, associated with transient shock distortion during vortex passage.
This feature remains confined to the shock and persists after the interaction.

Downstream of the shock, the vortex remains visible as a coherent structure, though it is deformed into a sheared core.
The solution remains stable throughout the interaction, with regularization localized to compressive regions and no visible spurious oscillations.

\subsubsection*{Supersonic flow interaction with curved solid boundaries}

Supersonic flow past a cylinder is used to examine the behavior of the regularization in the presence of boundary-induced compression and curved shock structures.
In contrast to the previous examples, where compression arises from prescribed shocks or shock interactions, the compression here is generated by the geometry itself.
The resulting bow shock is not aligned with the mesh and provides a test of robustness under curved boundaries.
The example examines whether the regularization remains localized to the shock while preserving stability and accuracy in the surrounding flow.

\begin{example}[Supersonic flow past a cylinder]
    \label{ex:supersonic_cylinder}
    The problem is posed on $\Omega=[0,4]\times[-1,1]$, with a circular cylinder of radius $R=0.25$ centered at $(0.6,0)$.
    The domain is discretized using a curved quadrilateral mesh constructed by mirroring an upper-half structured mesh about the centerline.
    The mesh includes a concentric annular layer around the cylinder and mild refinement near the outer slip walls.

    Inflow conditions are imposed at the left boundary using the free-stream state $U_{\mathrm{in}}$ with $(\rho_{\infty},u_{\infty},v_{\infty},P_{\infty})=(1.4,\,3.0,\,0.0,\,1.0)$.
    A supersonic outlet condition is applied at the right boundary.
    The top and bottom boundaries are treated as slip walls, and the cylinder surface is also treated as a slip wall.
    The interior initial condition is constructed from a smoothed start-up field near the cylinder walls.

    \noindent Let $r(x)=\sqrt{(x_1-0.6)^2+x_2^2}$, and, for $r(x)>0$, define the outward normal and tangent relative to the cylinder by
    \begin{equation}
        n(x)=\frac{1}{r(x)}\begin{pmatrix}x_1-0.6\\ x_2\end{pmatrix},
        \qquad
        t(x)=\begin{pmatrix}-n_2(x)\\ n_1(x)\end{pmatrix}.
    \end{equation}
    With free-stream velocity $\vct{u}_{\infty}=(u_{\infty},v_{\infty})^\top$, the normal and tangential free-stream components are
    \begin{equation}
        u_{n,\infty}(x)=\vct{u}_{\infty}\cdot n(x),
        \qquad
        u_{t,\infty}(x)=\vct{u}_{\infty}\cdot t(x).
    \end{equation}
    The start-up state uses the radial blending function
    \begin{equation}
        S(r)=
        \begin{cases}
            0, & r \le R,\\
            1-\exp\!\left(-\left(\dfrac{r-R}{\delta}\right)^2\right), & r>R,
        \end{cases}
        \qquad
        \delta=\max\{8h_{\min},0.08\},
    \end{equation}
    where $h_{\min}$ is the minimum mesh length scale.
    The initialized velocity field is then
    \begin{equation}
        \begin{aligned}
            u_n(x,0) &= S(r(x))\,u_{n,\infty}(x),\\
            u_t(x,0) &= u_{t,\infty}(x),\\
            \vct{u}(x,0) &= u_n(x,0)\,n(x)+u_t(x,0)\,t(x),
        \end{aligned}
    \end{equation}
    Density and pressure are initialized to the free-stream values, and the total energy is computed from the equation of state.

    Only the wall-normal component of velocity is damped near the cylinder, producing a smooth initial state consistent with the boundary conditions.
    As the flow evolves, a detached bow shock forms upstream of the cylinder.
    Geometric obstruction generates sustained compression and produces a curved shock that is not aligned with the computational mesh.
\end{example}

\begin{figure}[pos=htbp]
    \centering
    % \tikzsetnextfilename{supersonic_cylinder_2d_snapshots}
    % \input{figures/supersonic_cylinder_2d/supersonic_cylinder_2d_snapshots}
    \includegraphics{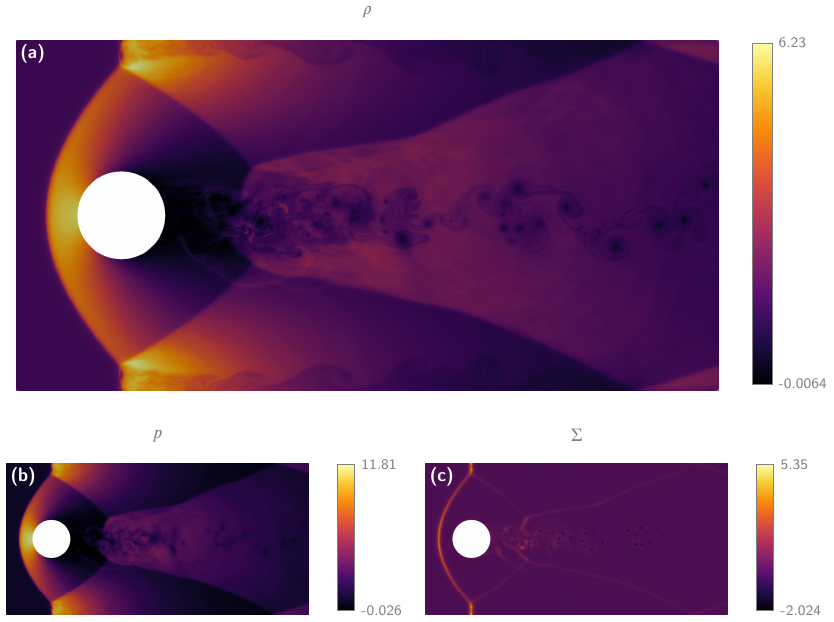}
    \caption{
        Two-dimensional supersonic cylinder problem (Example~\ref{ex:supersonic_cylinder}) for polynomial order $p=2$.
        Snapshots at $t=3.0$ are shown for (a) density $(\rho)$, (b) pressure $(P)$, and (c) entropic pressure $(\Sigma)$.
        Results are obtained using IGR--DG with three levels of refinement and $\alpha_{0}=1.0$.
    }
    \label{fig:supersonic_cylinder_2d_snapshots}
\end{figure}

As the solution evolves, the incoming supersonic flow interacts with the cylinder and forms a detached bow shock upstream of the body.
Representative snapshots of the density, pressure, and entropic pressure are shown in \cref{fig:supersonic_cylinder_2d_snapshots} at $t=3.0$.
The shock is curved and not aligned with the mesh, and $\Sigma$ activates sharply along the shock while remaining negligible in the free-stream region.
For the reported choice $\alpha_{0}=1.0$, the solution remains stable, although small oscillations are visible along the reflected shock fronts where the bow shock undergoes Mach reflection at the upper and lower walls.
The Mach reflection regions also exhibit additional fine-scale vortical features near the triple points.
Such secondary features generally require sufficient order and resolution to be resolved, and their persistence indicates that the regularization adds little dissipation away from compressive interfaces.

Downstream of the cylinder, the flow develops a low-pressure wake bounded by slip-line-type contact discontinuities emanating from the body.
The regularization does not directly target these features, since it is designed to activate in compressive regions rather than across contact discontinuities or other shear-dominated structures.
The Zhang--Shu positivity-preserving limiter~\citep{zhang2010positivity} is enabled in this case to maintain admissibility in the wake.
As a result, some smearing of fine-scale features is visible in the wake, and fine-scale vortical features are less sharply resolved.
Despite this, the overall flow remains stable.
The regularization remains localized to the bow shock and does not introduce additional dissipation in the free stream or along the solid boundary.
These results demonstrate robust behavior in the presence of curved geometry, boundary-induced compression, and non-aligned shock structures.

\newpage
\section{Comparison, conclusion, and outlook}
\subsection*{Comparison to prior work}

IGR introduces regularization at the level of the governing equations rather than through modifications of the discretization, distinguishing it from classical shock-capturing approaches for discontinuous Galerkin methods, including limiter-based methods and artificial viscosity.
Rather than stabilizing the numerical solution through reconstruction, limiting, or localized dissipation, IGR modifies the compressive dynamics through the entropic pressure obtained from an auxiliary elliptic equation.

Entropy-viscosity methods also regularize the governing equations through diffusion guided by entropy production~\citep{guermond2011entropy,cook2005hyperviscosity}.
IGR instead regularizes compressive dynamics through the entropic pressure without introducing viscous dissipation, producing a different balance between stabilization and accuracy.

The smooth convergence study shows that although IGR–DG reduces the order of accuracy before shock formation, it achieves the same order once shocks develop.
In post-shock flows, IGR–DG remains stable at higher polynomial orders while retaining higher-order structures.
At fixed degrees of freedom, increasing the polynomial order preserves these structures without increasing the error typically seen with methods that damp high-order features.

Compared with the characteristic TVB-limited DG formulation considered throughout this work, IGR--DG stabilizes shocks without shock-capturing limiters or artificial viscosity while retaining fine-scale flow features as the polynomial order is increased.

\subsection*{Conclusion and outlook}

In this work, we develop a discontinuous Galerkin discretization of the IGR-regularized Euler equations and evaluate its behavior in shock-dominated compressible flows.

Across a range of one and two-dimensional benchmark problems, the resulting IGR--DG method stabilizes shocks while the entropic pressure remains confined to regions of strong compression.
The regularized shock profiles remain narrow and sharply resolved, and the regularization remains localized even in the presence of shock interactions and boundary-induced compression.
Apart from localized oscillations observed along the reflected shock front in the supersonic-cylinder example, the method remains free of significant spurious oscillations in the benchmark problems considered.
Since IGR does not replace positivity-preserving methods, a positivity-preserving limiter is employed in the supersonic-cylinder example to maintain admissible states.

The numerical experiments demonstrate that fine-scale flow features are retained in challenging multiscale flows despite the presence of shocks.
This behavior is particularly evident in the Shu--Osher and high-frequency perturbation problems, where oscillatory structures remain well resolved away from shocks while maintaining stable shock resolution.

This work demonstrates that the proposed DG formulation effectively incorporates information geometric regularization into a high-order discontinuous Galerkin method.
The benchmark problems show that regularization through the entropic pressure provides stable shock resolution while retaining fine-scale flow features across a broad range of compressible flow problems.

\section*{Acknowledgments}
This work was supported by the Predictive Science Academic Alliance Program (PSAAP Award DE-NA0004261 - “The Center for Information Geometric Mechanics and Optimization
(CIGMO)”) managed by the NNSA (National
Nuclear Security Administration) Office of Advanced Simulation.
BE and FS gratefully acknowledge support from the Air Force Office of Scientific Research under award number FA9550-23-1-0668 (Information Geometric Regularization for Simulation and Optimization of Supersonic Flow).
FS was also supported by the Alfred P. Sloan Foundation via a Sloan Research Fellowship in Mathematics.

\bibliographystyle{cas-model2-names}
\bibliography{references}

\end{document}